\documentclass[11pt]{amsart}

\usepackage{amsmath}
\usepackage{amsfonts}
\usepackage{amssymb}
\usepackage{amscd}
\usepackage{color}
\input{epsf.sty}
\catcode `\@=11
\def\numberbysection{\@addtoreset{equation}{section}
         \renewcommand{\theequation}{\thesection.\arabic{equation}}}
\numberbysection
\def\subsubsection{\@startsection{subsubsection}{3}%
  \normalparindent{.5\linespacing\@plus.7\linespacing}{-.5em}%
  {\normalfont\bfseries}}
\newtheorem{introthm}{Theorem}
\newtheorem{thm}{Theorem}[section]
\newtheorem{lem}[thm]{Lemma}
\newtheorem{prop}[thm]{Proposition}
\newtheorem{cor}[thm]{Corollary}

\theoremstyle{definition}

\newtheorem{df}[thm]{Definition}

\newtheorem{rmk}[thm]{Remark}
\newtheorem{nota}[thm]{Notation}
\newtheorem{ex}[thm]{Example}

\newcommand{\aoa}[2]{a_{#1}\otimes \dots \otimes a_{#2}}
\newcommand{\bob}[2]{b_{#1}\otimes \dots \otimes b_{#2}}
\newcommand{\ata}[2]{a_{#1}\cdots a_{#2}}

\newcommand{\set}[2]{\{#1,\dots, #2\}}

\def\nn{\nonumber}

\def\In{In}

\def\Sn{\mathbb{S}_n}

\def\Snn{\mathbb{S}_{n+1}}

\def\Z{\mathbb{Z}}

\def\ra{\rightarrow}
\def\del{\partial}

\def\t{\tau}
\def\G{\Gamma}
\def\a{\alpha}
\def\b{\beta}
\def\g{\gamma}

\def\CHom{\mathcal{H}om}
\def\O{\mathcal{O}}

\def\pair{\la \;,\;\ra}

\def\del{\partial}

\def\a{\alpha}

\def\b{\beta}
\def\g{\gamma}

\def\t{\tau}

\def\nn{\nonumber}
\def\la{\langle}
\def\ra{\rangle}

\def\Sn{\mathbb{S}_n}

\def\Snn{\mathbb{S}_{n+1}}

\def\Hom{\mathrm{Hom}}

\def\comp{Comp}

\def\O{\mathcal{O}}
\def\P{\mathcal{P}}

\def\AP{\mathcal{P}^{\angle}}
\def\PA{\mathcal{P}^{\angle}}

\def\lab{\mathrm{lab}}

\def\arcgraphs{\overline{\mathcal{G}}}
\def\arcgraph{\overline{\mathcal{U}}}
\def\aarcgraphs{\arcgraphs^{\angle}}
\def\carcgraphs{\arcgraphs^{e}}
\def\qfarcgraphs{\arcgraphs_{\#}}
\def\Parcgraphs{\mathcal{P}\arcgraphs}
\def\Pcarcgraphs{\mathcal{P}\carcgraphs}
\def\Paarcgraphs{\mathcal{P}^{\angle}\arcgraphs}
\def\Pacarcgraphs{\mathcal{P}^{\angle}\carcgraphs}

\def\Darc{\mathcal{DA}rc}
\def\Arc{\mathcal{A}rc}
\def\Arcn{{\mathcal{A}rc_{\#}}}

\def\Arcno{{\mathcal{A}rc_{\#}^0}}

\def\Tree{\mathcal{T}ree}
\def\Lintree{\mathcal{L}Tree}

\def\OCDiarc{\mathcal{C}_o^*(\Diarc)}

\def\Gr{Gr}
\def\OC{\mathcal{C}_o^*}

\def\Poly{\mathcal{P}oly}
\def\Polydiag{\Poly_{\infty}}
\def\Ass{\mathcal{A}ssoc}

\def\Ana{\A^{\angle}}
\def\A{\mathcal{A}}

\def\Diarc{{\mathcal{A}rc^{i/o}}}

\def\DiA{\A^{i/o}}
\def\Diioa{\A^{i\leftrightarrow o}}
\def\Diioarc{\overline{\Arc}^{i\leftrightarrow o}}
\def\LDiioarc{\overline{\L\Arc}^{i\leftrightarrow o}}

\def\Diioarci{\overline{\Arc}^{i\leftrightarrow o}_1}

\def\iooarc{\mathcal{A}rc^{i \nleftrightarrow i}}
\def\ioarc{\mathcal{A}rc^{i\leftrightarrow o}}
\def\Anarc{\mathcal{A}rc^{\angle}}
\def\Anarcn{\mathcal{A}rc_{\#}^{\angle}}
\def\Anarcno{\mathcal{A}rc_{\#}^{0,\angle}}

\def\pmrib{\mathcal{R}ib^{part}}

\def\Rib{\mathcal{R}ib}

\def\PRib{\mathbb{P}\mathcal{R}ib}
\def\Anrib{\mathcal{R}ib^{\angle}}

\def\Loop{\mathcal{L}oop}

\def\mk{\mathrm{mk}}

\def\amark{mk^{\angle}}

\def\val{val}
\def\Ainf{A_{\infty}}

\def\mk{mk}
\def\io{i/o}

\def\In{In}
\def\Out{Out}

\def\lab{{Lab}}

\def\L{\mathcal{L}}

\def\Loop{\mathcal{L}oop}

\def\Arc{\mathcal{A}rc}
\def\Brace{\mathcal{B}race}

\def\Cyc{\curvearrowright}

\def\Rp{{\mathbb{R}_{>0}}}

\def\CC{CC} %cyclic chains
\def\CH{C} %hochschild chains
\def\HH{H} %hochschild homolgy
\def\HC{HC} %homology cylique

\def\TV{\overline{TV}}
\def\Shuff{\mathcal{CS}}
\def\Modshuff{\mathcal{CM}}
\def\Modshuffint{\mathcal{MCS}}

\def\Cyclic{\mathcal{C}yc}
\def\H{\mathcal{H}}

\def\Mngn{M_{g,n+1}^{1^{n+1}}}

\def\titem{\item[--]}

\newcommand\crl[1]{\vskip 5\lineskip \noindent {\bf Corollary. }{\em #1}
\vskip 5\lineskip}

\usepackage{amscd}

\begin{document}

\title[Moduli space actions on Hochschild co-chains II]
{Moduli space actions on the Hochschild Co-chains of a Frobenius
algebra II: Correlators}

\author
[Ralph M.\ Kaufmann]{Ralph M.\ Kaufmann}
\email{kaufmann@math.uconn.edu}

\address{University of Connecticut, Storrs CT 06269}

\begin{abstract}
This is the second of two papers in which we prove that a cell
model of the moduli space of curves with marked points and tangent
vectors at the marked points acts on the  Hochschild co--chains of
a Frobenius algebra. We also prove that a there is dg--PROP action
of a version of Sullivan Chord diagrams which acts on the
normalized Hochschild co-chains of a Frobenius algebra. These
actions lift to operadic correlation functions on the co--cycles.
 In particular,
the PROP action gives an action on the homology of a loop space of
a compact simply--connected manifold.

In this second part, we discretize the operadic and PROPic
structures of the first part. We also introduce the notion of
operadic correlation functions and use them in conjunction with
operadic maps from the cell level to the discretized objects to
define the relevant actions.

\end{abstract}

\maketitle

\tableofcontents

\section*{Introduction}
The subject of this sequence of two papers are actions on the
Hochschild complex of an associative or more restrictively a
Frobenius algebra induced by operadic structures on moduli spaces of
curves. Our approach is from the point of view of combinatorial
field theory which relies on a graph description of moduli space
that also allows for a natural compactification which was given by
Penner; see \cite{KS2,cost} for different approaches. One upshot of
our treatment is that the role of the Frobenius condition becomes
clear, thus allowing us to separate when this additional assumption
is needed and when one can get by with just an associative algebra.
Actions of a cell model of the open moduli space are expected on the
grounds of $D$--brane considerations \cite{KR,KLi1,KLi2}\footnote{A
more extensive discussion of these links is given in
\S\ref{conclusion}.}, while a subspace of moduli space given by
certain graphs is supposed to act by the considerations of string
topology \cite{CS,Vor,CJ,C,CG,Chat,Merk,C1,C2}. The archetype of
these actions was established with the proofs of Deligne's
conjecture \cite{Maxim,T,MS,MS2,MScosimp,Vor2,KS,MS2,BF,del} and its
generalizations to the $A_{\infty}$ \cite{KS,KSch} and to the cyclic
case \cite{cyclic}. The essential role of the topological operad in
this ``classical'' case was clearly assigned to the little discs and
the framed little discs operads. In the present setup the first
objective is to establish the existence of the topological and
cellular operads needed to make the above expectations into provable
statements. This task was completed in the first part \cite{hoch1}.
For instance one of the results of \cite{hoch1} is that there is a
rational operad structure on the chains of the moduli space $\Mngn$
of genus $g$ curves with $n$ punctures and a tangent vector at each
of these punctures which induces a chain level operad. This result
can be seen as a presentation of a combinatorial version of
conformal field theory (CFT) in terms of foliations \cite{KP}. As
explained in \cite{hoch1}, in our setup the operad structure for the
moduli space on the topological level  cannot be expected to be
strict, since it is only well defined almost everywhere. This is
captured by the notion ``rational''. Likewise, in \cite{hoch1}, we
showed that there are topological and chain level operads/PROPs for
the subspaces of Sullivan--Chord diagrams which are at the heart of
string-topology like operations. Here for the PROP structure we need
a weakening to a ``quasi''--PROP, which means that the associativity
only needs to hold up to homotopy. An important result of the first
part \cite{hoch1} is that these weaker structures nevertheless
induce the strict operad/PROP structures on the cell level.

In the present second part, we establish that the cell level
structures indeed act on the Hochschild co--chains. Our main tool is
the notion of operadic correlation functions, which should be
understood as a suitable definition of a dg--algebra $(A,d)$ over a
cyclic operad. Another way to phrase this is that these correlation
functions reflect  the fact that OPEs in physics are actually only
defined within correlators and only on BRST closed fields.

In particular, using the results of \cite{hoch1}, we prove the
following theorem announced in loc.\ cit.\:

\begin{introthm}
\label{introa}
 The moduli space $\Mngn$ of genus $g$ curves with
$n$ punctures and a tangent vector at each of these punctures
has the structure of a rational cyclic operad. This structure
induces a cyclic $dg$ operad structure on a cell model computing the
cohomology of $\Mngn$.

Furthermore the cell level operad operates on the Hochschild
co--chains of a Frobenius algebra. It also yields correlation
functions on the tensor algebra of the co-cycles of a differential
algebra $(A,d)$
 with a cyclically invariant trace $\int: A\to k$
which satisfies $\int da=0$ and whose induced pairing on
$H=H(A,d)$ turns $H$ into a Frobenius algebra.
\end{introthm}

As stipulated in \cite{hoch1} there is also a PROPic version of this
action involving (a partial compactification of) a subspace.  The
corresponding theorem pertaining to string topology type operations
is again proved in this second part.

\begin{introthm}{There is a rational topological quasi--PROP which
 is homotopic to a CW
complex whose cellular chains are isomorphic as a free Abelian
group to a certain type of Sullivan Chord diagrams. These chains
form a dg--PROP and hence induce this structure on the Chord
diagrams. Furthermore if $H$ is a Frobenius algebra there is a
PROPic action on the Hochschild co--chains of $H$ that is a
dg--action. This dg-action of a dg-PROP on the dg--algebra of
Hochschild co-chains naturally descends to an action of the
homology of the CW-complex on the Hochschild cohomology of a
Frobenius algebra.

Moreover for $(A,d,\int,H)$ as in Theorem \ref{introa} the action
on $H$ is induced by correlation functions on the tensor algebra
of $A$ that yield operadic correlation functions on the tensor
algebra of the co--cycles of $A$ for any $(A,d)$ as above.

Finally, the $BV$ operator, which is given by the action of the
sub--PROP equivalent to the framed little discs operad, acts as in
\cite{cyclic}. Thus the BV operator for the action on the
Hochschild cohomology of $H$ is identified with Connes operator
$B$ under the identification of the Hochschild cohomology of a
Frobenius algebra with its cyclic cohomology of $H$.}
\end{introthm}

The application to the homology of the loop space of a
simply connected manifold then comes as an immediate consequence
 using Jones' \cite{jones} cyclic description
of the free loop space.

\crl {When taking field coefficients, the above action gives a
$dg$--action of a $dg$--PROP of
 Sullivan Chord diagrams on the $E^1$--term of a spectral sequence
converging to $H_*(LM)$, that is the homology of the loop space of
a simply connected compact manifold, and hence induces operations
on the homology of the loop space.}

If we are not dealing with an algebra $A$, but only with a
$dg$--vector space $(V,d)$ that has a pairing $\pair$ which is
symmetric and  satisfies $\forall v,w\in V: \la dv, w\ra +\la
v,dw\ra =0$, such that $H=H(V,d)$ is finite and the induced pairing
on it is non--degenerate, there is a still an action.

\begin{introthm}
The operads and PROPs above also act on the tensor algebra $TV$ of
a triple $(V,d,\pair)$ as specified above and yield operadic
correlation functions for the co-cycles of $TV$.
\end{introthm}

This action is different from the algebra case of Theorems A and B,
making the result interesting in its own right. It seems to be that
there are two basic strategies to obtain correlators from graphs on
surfaces, a ``multiplicative'' action for an algebra and an
``additive'' action for a vector-space. One added feature of this
action ``descends'' to an action of the stabilized arc operad, which
forms a spectrum \cite{Ribbon}.

As mentioned above, the proof of these facts consists mainly of two
steps. First defining the respective topological objects and then
defining their actions. The first step is contained in \cite{hoch1}
and the second one is the content of this article, in which we
define the actions of the various objects. The definition of the
action itself again breaks down into two parts.

The first part of this paper is the very  definition of the
operations. In this aspect the paper is completely self--contained.
The approach we use is to first introduce discretized versions of
the topological operads and PROPs of \cite{hoch1} and then to let
these operate via correlation functions. In this completely cyclic
setting it is more natural to define multilinear maps to $k$ rather
than maps in $\Hom(V^{\otimes n},V)$. The problem with this approach
is that maps to $k$ are not easily composed which is why we
introduce the notion of operadic correlation functions. The approach
of using correlation functions also mixes well with the ideas of
physics where these objects are fundamental. Taking a physics
perspective, anything which does not change a correlation function
is not physical; or in other word OPEs live only inside correlators.
In one of our main examples, namely that of a quasi--Frobenius
algebra, this means that we can lift the constructions from the
cohomology to the co--cycle level.

The second part is to show that this action has the desired operadic
or PROPic properties. For this we use the operadic correlation
functions to get the results on the discretized level. The last
step, which is the one that requires the results of \cite{hoch1} is
to relate the operadic/PROPic structure of the discretized/graph
level to the cell level operations of the chain model operads whose
theory is developed in \cite{hoch1}. The relevant facts are reviewed
in the first paragraph. The language we use is that of arc--graphs
on surfaces. This ties in with the description of the arc--operad
$\Arc$ of \cite{KLP}. In the special case of quasi--filling arcs,
that is the subspace $\Arcno$ of $\Arc$ which corresponds to the
moduli space $\Mngn$, there are actually two formalisms which one
can use: The arc graphs and their dual ribbon graphs. We write out
the details in both of these pictures, so that the reader more
familiar with ribbon graphs can more easily understand the
constructions. They are however different from the usual known
constructions and they do not generalize to the boundary, that is to
the more general, non--quasi--filling case, which is needed to
define the String--Topology type operations.

When dealing with the operadic/PROPic properties, one has to be very
careful about the operations on the side of the endomorphism operad
$\CHom$ of the relevant linear spaces. This is a subtlety which is
know from Deligne's conjecture. When being precise about the signs,
one actually does not prove  that one has an algebra over the
relevant operad, but rather an operadic morphism to the operad
$\Brace$ which is formed by subspaces of the endomorphism operad,
but has different sign rules and hence a twisted operad structure.
When dealing with our actions, a similar situation arises which is
slightly more complicated. We again obtain an operadic morphism to
an operad which is formed by subspaces of the endomorphism operad.
These spaces have a grading and the induced sub-operad structure
respects the associated filtration, but not the grading. Projecting
to the associated graded operad structure  and  correcting the sign
according to the grading, we obtain operadic actions as operadic
morphisms to these ``twisted'' endomorphism sub--operads.

The paper is organized as follows: In the first paragraph, we review
the necessary facts we need from \cite{hoch1} and then define the
``discretized'' versions of the spaces we will consider. These
``discretized'' versions are free Abelian groups of graphs on
surfaces, so--called partitioned arc graphs. We then define operad
and PROP structures on these graphs and go on to show
 that partitioning an arc--graph, which by \cite{hoch1}
(also see below) can be thought of as indexing a cell of an operadic
cell complex, is an operadic morphism. This is actually quite
subtle, since different types of graphs require different types of
cell operads. The principal choices are filtered or graded versions.
There is an intricate interplay between the discrete data associated
to the graphs and the geometry they realize. In \S3 we define the
notion of operadic correlation functions and give several examples.
The next paragraph \S4 is dedicated to defining correlation
functions, aka.\ correlators, for an all encompassing class of
graphs, the angle marked partitioned arc graphs. These correlators
are actually defined on the tensor algebra of a quasi--Frobenius
algebra. In \S5 we show that the correlators become operadic in
several different settings. Notably for $\Diioarc$ and for $\Arcno$,
see \cite{hoch1} and \S\ref{spacelistpar} below. In the latter case,
we have to be careful about the operadic structure of the spaces the
operations take values in. The relevant subset of the $\CHom$ operad
is graded and hence filtered. As mentioned above, the correlators
 define an operadic map to the associated graded of
this filtration. {\em A priori} the operadic compositions in the
$\CHom$ operad and the discretized graphs only agree up to lower
order terms in the associated filtration. {\em A posteriori} these
terms agree for the action of $\Diioarc$. The last paragraph \S6
contains concluding remarks about the link to $D$--branes and future
research directions.

\section*{Acknowledgments} We would like to thank the Max--Planck--Institute
for Mathematics where this work was started, a good portion of
it  was written in the summer of 2005 and the finishing
touches were put on in the summer of 2006.  The two papers received their
final form at the MSRI, which we would like to thank for its
hospitality in May 2006. It is a pleasure to thank
 Bob Penner, Ralph Cohen, Jim McClure, Dev Sinha and Craig Westerland for
 discussions
on various details during various stages of this project.

\section*{Conventions}
We fix $k$ to be a field of arbitrary characteristic.
Also in this
part of the paper we always assume that the number of punctures is
zero. That is $s=0$ for all operads and suboperads.

\section{Discretizing the Arc operad and its cousins}
\subsection{Brief review}
\label{review} Without going into the details, which are contained
in \cite{hoch1}, we wish to point out the basic definitions of the
graphs underlying the various versions and generalizations of the
$\Arc$ operad. On a proper subset of $\Arc$, the quasi-filling
arc--families $\Arcno$ there are two pictures, one in terms of arc
graphs and one in terms of the dual ribbon graphs. Although this
subset is not big enough, even for our purposes ---for instance to define the
string topology type operations--- we include both pictures, since
ribbon graphs are commonly used to describe moduli spaces and are
hence predominant in the literature.

\subsubsection{Graphs}
 A graph  is a tuple
$\G=(V(\G),F(\G),\imath_{\G},\del_{G})$ where $V(\G)$ is a set
whose elements are called the vertices, $F(\G)$ is a set whose
elements are called the flags or ``half edges'', $\imath_{\G}:
F(\G)\to F\G)$ pairs the ``half edges'' to edges and
  $\del_{\G}:F(\G)\to V(\G)$ gives the vertex of a flag. An edge in
this setting is an orbit of $\imath$ that is a set of flags
$\{f,\imath(f)\}$. An oriented edge is a pair of flags
$(f,\imath(f))$. The set of flags incident to a vertex $v$ is
called $F_{v}(\G)$.

Recall that a ribbon graph is a graph with a cyclic order of each
of the sets of flags incident to a fixed vertex. Such a ribbon
graph has natural bijections $\Cyc_v:F_v\rightarrow F_v$ where
$\Cyc_v(f)$ is the next flag in the cyclic order. Since $F=\amalg
F_v$ one obtains a map $\Cyc:F\rightarrow F$. The orbits of the
map $N:=\Cyc \circ \imath$ are called the cycles or the boundaries
of the graph. These sets have the induced cyclic order. Due to the
cyclic order a ribbon graph also can be ``fattened'' to a surface
with boundary, by realizing the graph as a CW complex and then
thickening the edges to bands. In this fashion one obtains a
surface whose boundary components correspond to the cycles. The
genus of such a graph is given by the genus of this surfaces.
Explicitly, $2-2g(\G)=|V(\G)|-|E(\G)|+\#\text{cycles}$. An $S$
marking of a ribbon graph $\G$ is a bijection
 $\{cycles (\G)\}\to S$.

An angle is a pair of flags $(f,\Cyc(f))$, we denote the set of
angles by $\angle_{\G}$. It is clear that $f\mapsto (f,\Cyc(f))$
yields a bijection between $F_{\G}$  and  $\angle_{\G}$ . An angle
marking by a set $T$ is a map $\amark:\angle_{\G}\to T$. We will
call a (not--necessarily connected) ribbon graph with an angle
marking by $\Z/2\Z$ simply an angle marked ribbon graph.

\subsubsection{Arc graphs}

Fix an  oriented surface $F=F_{g,r}^s$ of genus $g$, with $s$
punctures and  $r$ labelled boundary components that each contain
one marked point. We usually label the boundaries from $0$ to $r-1$.
An arc graph on $F$ is a class of graphs on $F$. It can be thought
of as the orbit of a graph whose vertices coincide with the marked
points on the boundary under the action of the pure mapping class
group which fixes the marked points on the boundary and the
punctures pointwise. The edges of the graphs which comprise the
orbit are embedded arcs considered up to homotopy. We frequently
call these edges ``arcs''. There are certain conditions on the
graphs whose orbits we consider
\begin{itemize}
\item[i)]at least one arc. \item[ii)]  no parallel arcs (by
homotopy fixing the endpoints) \item[iii)]  no arcs parallel to a
boundary component (again using a homotopy fixing the endpoint)
\end{itemize}
We wish to point out that the underlying abstract graph of an arc
graph which is given by the collection of vertices and edges
together with their incidence relations and the set of complementary regions
of the arcs is invariant (or better equivariant) under the
action of the mapping class group. When we depict arc graphs, we
choose a particular representative. The class of all of the graphs
is called $\arcgraphs$. Since a vertex of the arc graph
corresponds to a boundary component, the vertices are labelled. We
frequently write $v_i$ for the vertex labelled by $i$, that is the
unique vertex lying on the boundary $i$ of the surface. Strictly
speaking an arc graph is a triple $(F,\G,\overline{[i]})$ of a
surface $F$, a ribbon graph $\G$ and an orbit of a homotopy class
of embeddings of the graph into the surface $\overline{[i]}$. The
full details are contained in \cite{hoch1}.

An arc graph is called {\em exhaustive} if there are no vertices
with valence $0$ and {\em quasi-filling} if the complementary
regions of the arcs are at most once-punctured polygons. The former
class of arcs is called $\carcgraphs$ and the latter class of arc
graphs is called $\qfarcgraphs$. An arc--graph becomes a (possibly
disconnected) ribbon graph by using the orientation on the surface.
It moreover even has a linear order of all the flags at a vertex due
to the induced
 orientation on the boundaries of the surface and hence a total order
 on all flags, by first enumerating the flags according to their labelled boundary
 components and then according to their linear order at that component.

An arc graph is called twisted at the boundary $i$ if the first and last
arc incident to $i$ are homotopic in $F$, when one allows homotopies that
move the endpoint on the boundary component $i$.

To each  quasi-filling arc graph there is  dual graph which is a
marked ribbon graph.  To simplify say the graph is on a surface
$F_{g,r}^0$, viz.\ no internal punctures.   The dual graph is then
defined as follows: Choose a representative $\g$, decompose the
surface into the complementary regions that is the components of
$F\setminus \g$. Now associate a vertex to each complementary
region and an edge to each arc. The edge is fixed to connect the
vertices (or vertex) representing the regions on the two sides of
the arc. This graph is again a ribbon graph, by using the
orientation of the arc graph to give the arcs bordering a
polygonal complementary region a cyclic structure. Each cycle of
this dual ribbon graph corresponds to a boundary component of the
surface and hence has a {\em linear} order. That is for each cycle
there is a distinguished flag which is the first in this cycle.

Vice--versa, we can ``fatten'' the dual ribbon graph to a surface
and consider the graph as the spine of this surface.
 Applying a dual graph construction
in this setting produces an inverse to the construction of the dual
graph of an  arc graph (see \cite{hoch1} for full details). This
explains the terminology ``dual graph''; the case with no punctures
is all we will use in the following considerations. The case with
punctures is treated in \cite{hoch1}. The map $\Loop$ of \cite{KLP}
gives a generalization of the dual graph to the non-quasi-filling
case.

In an arc graph, not all the angles are on an equal footing. The
last and first flag at a vertex form a distinguished angle which
is called the outside angle at that vertex.  All angles beside the
outside angles are called inner angles.

We will also
 consider arc graphs in which the set of boundaries of the surface
is partitioned into the sets $\In$ and $\Out$.
 This partitioning is encoded in a map $\io:V(\G)\to \Z/2Z$,
where the value $1$ stands for ``in'' and $0$ stands for ``out''.
Recall that the set of vertices of the arc graph can naturally be
identified with the boundary components of the surface.

In the dual ribbon graph picture, we accordingly
have a labelling of the {\em cycles}
by $\Z/2\Z$ indicating ``in'' and ``out''.

\subsubsection{Spaces of graphs}
\label{spacelistpar} We obtain the space of a given class of
graphs, by looking at the set of projective metrics, that is
equivalence classes of maps $w_{\G}:E(\G)\rightarrow {\mathbb
R}_{>0}$ under the action of $\Rp$ by a global re-scaling; that is
the action given by $\lambda \in \Rp: (\lambda w)(e)=\lambda
w(e)$. The set of all graphs of a given class with projective
metric basically gets a topology by identifying the limit in which
$w(e) \to 0$ for some edge with the graph in which $e$ is deleted
(see \cite{hoch1} for details). We usually call elements of these
spaces projectively weighted arc-families in keeping with
\cite{KLP} and the work of Penner.

The most important spaces are:

\begin{itemize}
\titem $\Arcno$ the space defined by quasi--filling graphs. This
space is isomorphic to $\Mngn$, the moduli of curves of genus $g$
with $n$ marked points and one tangent vector at each of these
points.

\titem $\Diioarc$ the space of arc graphs with a projective
metric, together with a partitioning $\io$ into $\In$ and $\Out$
which satisfy the conditions (1) only arcs between ``in'' and
``out'' and (2) each ``in'' boundary vertex has valence at least
1. This space plays the role of Sullivan Chord diagrams.
\end{itemize}
A reference list of the  spaces  that will make an appearance are:

\begin{itemize}
\titem $\A$ the space of all arc graphs with a projective metric.
This is a CW complex whose cells are indexed by the arc graphs.
\titem $\Arc$ the sub--space
of all exhaustive arc graphs with a projective
metric.

\titem $\A^{\angle}$   the space of all elements of $\A$
 with an additional angle marking.
\titem  $\Anarc$ the space of all the exhaustive arc graphs together
with a projective metric and an angle marking by $\Z/2\Z$. We will
consider $\Arc$ as a subspace of $\Anarc$ by choosing the constant
marking $\amark\equiv 1$.

\titem  $\DiA$ the space of arc--graphs with projective metric
which have an additional marking $\io$ of the boundaries
distinguishing inputs and outputs. We will consider this space
again as a subspace of $\A^{\angle}$, by marking all outside
angles and all inner angles of the $\In$ boundaries by $1$ while
marking all marking inner angles of the $\Out$ boundaries by $0$.

 \titem
$\Diioa$ the sub-space of $\DiA$ which is comprised of the arc
graphs with a projective metric,
 that additionally
satisfy the condition that there are only arcs between $\In$
boundaries and the $\Out$ boundaries.

\titem $\ioarc$ the subspace whose underlying arc graphs are
exhaustive and all of whose arcs only run from $\In$ boundaries to
$\Out$ boundaries.

\titem $\Diioarc$ the subspace of $\Diioa$ whose underlying arc
graphs hit all the $\In$ boundaries.
\end{itemize}

The spaces above  naturally come as disjoint union over the number
of boundary components, which we usually think of  as labelled by
$\set{0}{n}$. In the case of $\DiA$ we first label the boundaries
and $\In$ and $\Out$ and then label these boundaries separately,
say, by $\set{1}{n}$ and $\set{1}{m}$. There are natural actions
of the permutation groups on these labels. In \cite{hoch1} we
showed that essentially that $\Ana$ and its subspaces are operads
and that $\Diioarc$ is a quasi--PROP. Actually some of these
spaces, notably $\Arcno$ are only rational operads, viz.\ defined
on an open dense set. The full details are quite elaborate and
make up the bulk of \cite{hoch1}.

\subsubsection{Operads/PROPs of arc graphs}

Each of these spaces has an associated graph--complex--cell--model
given by considering the free Abelian groups generated the
underlying graphs.
 The natural differential is given by restricting
the differential of the CW complex $\A$. This differential applied
to a graph is  the sum of the arc--graphs obtained by removing one
arc with the appropriate sign. The differential for the open cells
of a subspace is defined to be the sum over only those graphs which
correspond to elements in the subspace. In other words these
complexes are the relative complexes of the subspaces in $\A$. These
complexes inherit   algebraic operad and/or PROP structures by
treating the graphs as ``open cells''. Here an open cell is
constituted by elements corresponding to the possible projective
metrics of a fixed arc graph and we label the generator
corresponding to an open cell by the respective graph. Then there
are induced gluing operations from the topological level on the
``open cells'', see \cite{hoch1} for details. These cells are graded
by their dimension, which is the number of edges of the graph minus
one. We note for later, that the number of edges is also twice the
number of flags which coincides with the number of angles. The
gluing operations respect the filtration induced by the grading and
accordingly we obtain two versions of cell operads on the
graph--complex. The first is the induced structure on the ``open
cells'' which we denote by $\OC( . )$ and second one is the one
induced by the first structure on the associated graded of the
filtration by dimension. The latter is again of course additively
isomorphic to the former and both are isomorphic to the graph
complex. The operations differ however. To make this distinction
clear we denote the graph--complexes with the operations
corresponding to the associated graded by $\Gr\OC( . )$. On the cell
level possibly after passing to the associated graded, we always
obtain the honest structure, that is not the up to homotopy or a
rational version. Most importantly:

\begin{itemize}
\titem
 The associated graded cell complex
$\Gr\OC(\Arcno)$ is a $dg$ operad and $\Gr\OC(\Arcno)$ calculates
the cohomology of $\Mngn$. By using the angle marking $\amark
\equiv 1$ the operad $\Gr\OC(\Arcno)$ embeds into
$\Gr\OC(\Anarc)$.
\titem The cell complex $\OC(\Diioarc)$ is a PROP and
$\Gr\OC(\Diioarc)$ isomorphic to the cellular chains
$CC_*(\Diioarci)$ of a CW-complex $\Diioarci$ and these chains
have the structure of a $dg$--PROP. This PROP can be thought of as
the PROP of Sullivan--Chord diagrams.
\end{itemize}
The details are quite involved, and they are carefully written out
in \cite{hoch1}.

\subsection{Partitioned Ribbon graphs}
\subsubsection{Inserting Points into edges}
To define the operations on the Hochschild co-chains, we will
systematically deal with unstable graphs i.e.\ graphs which have
vertices of valence two. For this we will need the operation of
inserting vertices valence 2 into edges and also the reverse
operation of removing them.

First, we recall the notion of a marked ribbon graph, viz.\ the
type of stable graphs we consider.
\begin{df}
 A {\em marked ribbon graph} is a
ribbon graph together with a map $\mk:\{cycles\} \rightarrow
F_{\Gamma}$ satisfying the conditions
\begin{itemize}
\item[i)] For every cycle $c$ the directed edge $\mk(c)$ belongs
to the cycle.

\item[ii)] All vertices of valence two are in the image of $\mk$,
that is $\forall v,\val(v)=2$ implies  $v\in Im(\del\circ\mk)$.
\end{itemize}
\end{df}
We called the set and the Abelian group generated by these graphs
$\Rib$. It is naturally the disjoint union over the graphs
$\Rib(n)$ which have $n$ cycles. We showed in \cite{hoch1} that if
we think of $\Rib$ as labelled graphs, they form an operad by
inducing said structure via the isomorphism $\Rib\cong
\Gr\OC(\Arcno)$ induced by the dual graph construction.

\begin{df}
Given a graph $\G$ and an edge $e=\{f_1,f_2\} \in \G$ we let
$\G(e')$ be the graph whose vertices are $V_{\G}\amalg \{v_2\}$ with
flags $F_{\G}\amalg \{n_1,n_2\}$ with $\del(n_{i})=v_2$ and
$\imath(f_i)=n_i$. We say $\G(e')$ is obtained from $\G$ by
inserting a vertex into $e$. Notice that this insertion does no
disturb the cycles, that is there is a canonical identification of
the cycles before and after the insertions. If $\G$ has a marking on
its cycles this marking will simply be retained.

Vice-versa, if $v_2$ is a vertex of valence 2 with flags $n_1,n_2$
with $\imath(n_i)=f_i$ then we let $\G/v_2$ be the graph whose
vertices are $V_{\G}\setminus \{v_2\}$, whose flags are
$F_{\G}\setminus \{n_1,n_2\}$ and whose new relation for $\imath$ is
$\imath_{\G/v_2}(f_1)=f_2$. In case $n_i=\mk(c)$ for the cycle $c$
it lies on, we set $\mk_{\G\setminus
v_2}=(n\circ\imath)^{-1}_{\G}(n_i)$.

We write $\G'\triangleright \G$ if $\G$ is obtained from $\G$ by
repeatedly inserting vertices, i.e.\ if there is a sequence
$\G^{0},\dots, \G^{n}$, $\G^{i}=\G^{\prime i-1}(e)$ for some $e\in
E_{\G^{i-1}}$, and $\G'=\G^{n},\G=\G^{0}$.
\end{df}

\begin{df}
A {\em partitioned marked ribbon graph}
 is a
ribbon graph together with a map $\mk:\{cycles\} \rightarrow
F_{\Gamma}$ which satisfies the condition that for every cycle $c$
the flag $\mk(c)$ belongs to that cycle.

We let $V_2(\Gamma)=\{v\in V_{\G}, \val(v)=2\}$ be the vertices of
valence two and set $V_{part}=V_2\setminus Im(\del\circ\mk)$ to be
the partitioning vertices.
\end{df}

\begin{nota}
\label{shorthand}
Let $\pmrib$ be the set of all partitioned marked ribbon graphs.
To avoid cluttered notation, as we have done in \cite{hoch1}, we
abuse notation and denote by $\pmrib$ the set of graphs, the
Abelian group generated by it, as well as the collection
$\{\pmrib(n)\}$ where $\pmrib(n)$ is the Abelian group generated
by the subset of $\pmrib$ of graphs which have $n+1$ cycles that
are labelled by $\set{0}{n}$ together with the $\Snn$ action
permuting these labels. The various meanings will always be clear
from the context.
\end{nota}

For a marked ribbon graph $\G$ we will consider $\P:\Rib
\rightarrow \pmrib$
\begin{equation}
\P(\G)=\sum_{\G' \triangleright \G} \pm\G'
\end{equation}
and call it the partitioning of $\G$.
 The right hand side is
infinite, but it is graded by the number of partitioning vertices
$|V_{part}|$. The fastidious reader hence may take the $\P(\G)$ to
lie in $\pmrib[[t]]$, where $t$ is a variable whose power
corresponds the number of partitioning vertices. The sign is
explained in \S\ref{signpar}.

\begin{rmk}
This construction can be seen as  the generalization of the foliage
operator of \cite{del,cyclic} from the setting of treelike ribbon
graphs to the setting of general ribbon graphs.
\end{rmk}

\subsection{Discretizing the $\Arc$ operad}
\label{discpar} Given a tuple $\a= (F,\Gamma,\overline{[i]})$, we
will consider a series of embedded graphs which are obtained by
drawing parallel arcs.

\begin{nota}
Let $p=(n_1,\dots,n_k)$ be an ordered partition of $n\in
\mathbb{N}$ with each $n_i>0$ and denote the set of all these
partitions by $P(n,k)$.
\end{nota}

Let $\a\in\carcgraphs$.  Recall that all edges $E(\a)=E_{\G(\a)}$
are linearly
 ordered by enumerating the flags in the following order, first according to their boundary
 and then according to the linear order at that boundary component
 induced by the orientation of the surface.
 Let $k=|E(\a)|$ and $p\in P(n,k)$. We define $\a^p$
to be the  embedded graph obtained from$(F,\Gamma,\overline{[i]})$ by inserting $n_i-1$ parallel edges to
the $i$-th edge $e_i$ and embedding them parallel to $e_i$. We
call the result of this operation a partitioned arc graph and
denote the set of these by $\Pcarcgraphs$. And again we use the
conventions explained in \S\ref{shorthand}.

\subsubsection{Drawing arc--graphs}
\label{drawpar} 
\begin{figure}
\epsfbox{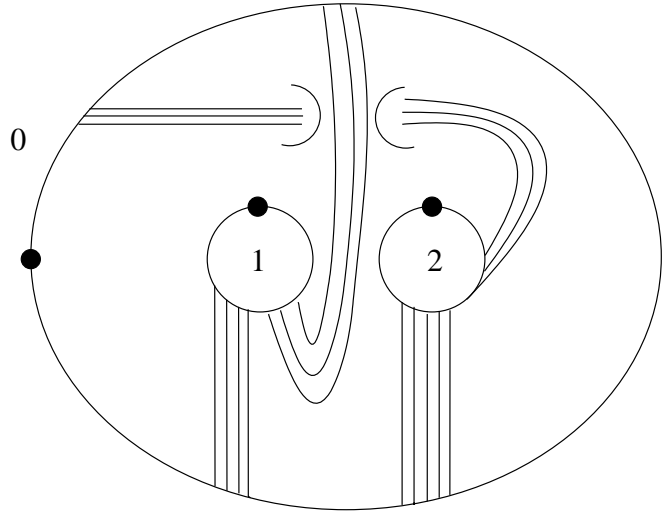}
\caption{\label{g2}
A  partitioned arc graph for the partition $(4,5,3,3)$}
\end{figure}

An example of such a graph is given in Figure
\ref{g2}. In this figure and in all other figures, we have taken the
liberty to depict the arc--graphs in a more suitable way. By
definition, all edges of a graph are incident to the vertices, which
would clutter the pictures. We therefore move the endpoints along
the boundary component that they are incident slightly apart. This
is done in such a fashion, that they (1) all are distinct and
distinct from the original vertex, and (2) their linear order along
the boundary component starting at the original vertex coincides
with their original linear order. We could have even defined the arc
graphs in this manner; see e.g.\ \cite{KLP} for all topologically
equivalent ways to define the space $\A$. We will therefore
henceforth use both pictures: the one with the edges apart as in
Figure \ref{g2} and the true graph picture, where all the edges are
incident to their vertices. The advantage of the latter lies in the
more direct definition and the advantage of the former is twofold,
first one obtains nicer pictures and secondly, the boundaries of the
complementary regions are $2k$--gons whose sides alternatingly
correspond to boundary components and arcs. This last observation
will make the definition of the action of these graphs on the
Hochschild co--chains more transparent.

We define
\begin{equation}
\label{Peq} \P(\a)=\sum_{n\geq k}\sum_{p\in P(n,k)}\pm\a^p
\end{equation}
As always when dealing with this type of object, one can use,
 a formal variable $t$ to keep track of the total number of
edges and consider the expression of equation (\ref{Peq}) as a
formal power series in $\PRib[[t]]$. The sign is explained in
\S\ref{signpar}.

\subsubsection{The underlying arc graph}
Given a partitioned arc graph $\g\in \P\arcgraphs$ it is possible to
recover  a unique $\a\in \arcgraphs$ such that $\g$ is a summand on
$\P(\a)$. Namely, for  $\g\in \P\arcgraphs$, we define the
underlying arc graph $\arcgraph(\g)$ to be the graph obtained by
gathering all parallel edges into one edge. Formally we introduce an
equivalence relation on $E(\g)$ where $e_1\sim e_2$ if $e_1$ and
$e_2$ are parallel ---as arcs in the surface via a homotopy fixing
the endpoints-- and then take the quotient by this relation. This is
well defined since if $e_1\sim e_2$ then $\del(e_1)=\del(e_2)$.

We call a partitioned arc graph {\em quasi--filling} if the
complementary regions are either polygons or once--punctured
polygons. Notice that a partitioned arc graph is quasi--filling if
and only if its underlying arc graph is quasi--filling.

\subsubsection{Grading}
Let $\Pcarcgraphs(n)^k$ be the set of those graphs whose underlying
arc graph has  $k+1$--edges which ``live'' on surfaces with $n+1$
 boundaries labelled from $\{0,\dots ,n\}$. We let
$\Pcarcgraphs(n):=\bigoplus_k \Pcarcgraphs(n)^k$. This space is
then graded by $k$ and hence filtered by the elements of degree
$\leq k$.

\begin{rmk}
We can view $\Pcarcgraphs$ as a discretized version of $\Arc$ in the
following two ways: $\P(\a)$ can be thought of as either (A) as a
sampling by the numbers $\frac{i}{k}$ of the boundary considered as
the interval defined by the window on the boundary (see
\cite{KLP,KP} for the formalism of windowed surfaces) or
equivalently (B) as a cosimplicial realization of the arc graph
$\a$. We wish to pursue the latter point of view elsewhere (see also
the comments  \S\ref{conclusion}  below).
\end{rmk}

\subsubsection{Signs} \label{signpar}
 In this paragraph, we wish to fix our sign conventions once and for all.
 As in \cite{del}, we do this by using tensor products indexed by sets in the spirit of
 \cite{KS} which makes all signs
 completely natural. We will henceforth not bother with them
 again. For other different explicit sign fixing schemes for operations
 of cell operads, we refer to
 \cite{del}.

First notice that for any arc graph partitioned or not, there is a
natural linear order on all the flags and hence on all the edges.
So we can use
 these linear orders to fix the signs.
In general one can do this quite nicely by considering
 the tensor product over $\mathbb Z$ of
the generator given by the graphs $\a$ with copies of a
 ``line of degree $1$'' that is a freely generated Abelian
group generated by an element of degree $1$. Thus we replace $\a$
with the expressions $\a\otimes L^{\otimes S}$ where $S$ is an
indexing set. In the cyclic operad setting we will use $S=E'(\G)$
that is the set of edges without the last edge in the linear order.
This gives a universal way to fix the signs. It also
 assigns the correct degree to $\a$ if $\a$ is thought
of as an element of the various cellular chain complexes
 introduced in \cite{hoch1}.
For the expression above this results in the sign obtained from the
shuffle $L^{\otimes E'(\a)}\otimes L^{\otimes n}$ to $L^{\otimes
E'(\a^p)}$ for the summand indexed by $p\in P(n,k)$.

 In the PROP setting the natural indexing set will be
$S=\angle^{in}_{\rm inner}(\G)$ the set of inner angles on the
inputs. This again corresponds to the dimension of the cells when we
consider the CW complex $\Diioarc$.

Of course one of the two sign conventions can be obtained from the
other by shifting the complexes. Alternatively, one could grade by
shifting by the dimension of the corresponding spaces to get rid of
the signs for top dimensional cells, as was done to obtain the Hopf
algebra of Connes and Kreimer in \cite{del}. Or, one could shift by
the number of all the boundaries. In this vein, we can consider the
use of $E'$ as the shift from the grading of $\Darc$ by a line
associated to the operadically distinguished boundary $0$.

%A slight variation of this theme is to use lines of degree one also
%for the boundaries. For instance one can assign one ``line of degree
%$1$'' to each edge and a line of degree  $-1$ to each $In$ boundary
%and a line of degree $1$ to each $Out$ boundary. This grading will
%then be the grading for the PROP. One can think of one equation

\subsection{Compatibility of the two constructions in the quasi-filling case}

Recall that for a quasi-filling arc graph $\G$ we defined a dual
ribbon graph $\hat \G$ in \cite{hoch1} as outlined in
\S\ref{review}. This construction easily generalizes to
quasi-filling $\G\in \Pcarcgraphs$. Moreover, since we are simply
inserting parallel edges, on can see that $\G \in \Pcarcgraphs$ is
quasi-filling if and only if $\G$ is a summand of $\P(\G')$ for
some quasi-filling arc graph $\G'$.

\begin{lem}
Denoting the dual ribbon graph of a quasi--filling
(partitioned) arc  graph $\a$
by $\hat \G (\a)$ and extending this construction linearly to give
a map of the respective Abelian groups, the following equality holds.
$$\hat \G(\P(\a))=\P(\hat \G(\a))$$
\end{lem}

\begin{proof}
 The insertion of parallel
edges corresponds to adding rectangles into the set of
complementary regions and this corresponds to inserting vertices
into the edges of the dual ribbon graph.
\end{proof}

\begin{figure}
\epsfxsize = \textwidth
\epsfbox{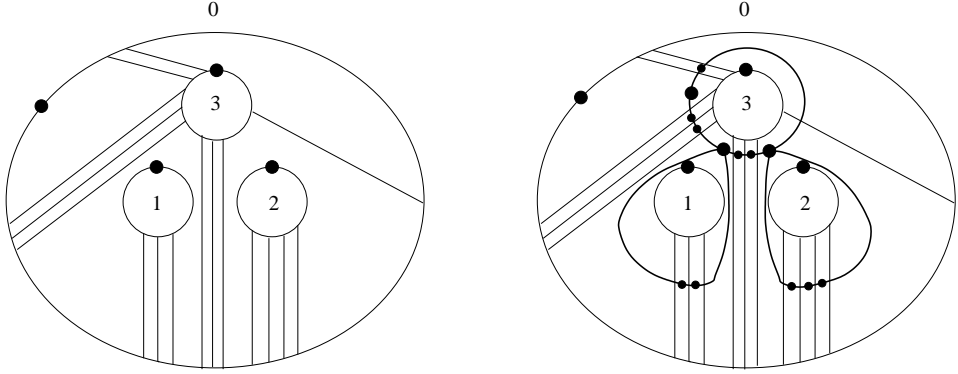}
\caption{\label{g2new}
A quasi--filling partitioned arc graph and it dual partitioned marked Ribbon graph}
\end{figure}

For an example of a partitioned graph and its dual see Figure
\ref{g2new}.

\subsubsection{The operad structure}
\label{partgluingpar} The set of Abelian groups generated by the
$\P\arcgraphs(n)$ has a natural cyclic operad structure, where the
$\Snn$ action is given by permutation of the boundary labels.

Given two elements $\a\in \P\arcgraphs(m)$ and $\b \in
\P\arcgraphs(n)$, we will call them {\em  matched} at the
boundaries $i$ and $0$, if the number of arcs incident to $\a$ at
the boundary $i$ is equal to the number of arcs incident to $\b$
at the boundary $0$. In case $\a$ and $\b$ are matched at the
boundaries $i$ and $0$, we define $\a\circ_i\b$ to be the graph
obtained by gluing the edges of $\a$ incident to $i$ to those of
$\b$ incident to $0$ matching them according to their linear
order. As usual this means that we reverse the linear order at
$i$, that is read off last to first, and match this order with the
linear order on $0$. If there is no matching at the boundaries, we
set $\a\circ_i\b=0$.
There are two final steps,

\begin{enumerate}
\item if both $\a$ and $\b$ are twisted at their boundaries $0$
and $i$, we set their composition to zero, and

\item if there are any closed loops, that is embedded arcs that do
not touch any boundary, as edges, we set the product to zero.
\end{enumerate}
If we omit step (1) we also call the gluing the algebraic gluing
and contrast call the gluing with both steps (1) and (2) the
geometric gluing, see \S\ref{naturalrmk} for comments.

It is clear that this gives an operad structure and that the subset
$\Pcarcgraphs(n)\subset \Parcgraphs(n)$ of exhaustive graphs is a
suboperad.

\begin{lem}
\label{Plem} The map $\P$ is an operadic morphism; that is we have
the following formula for the compatibility between the partitioning
of arc families and the operad compositions. For $\a,\b \in
\arcgraphs$, then
\begin{equation}
\P(\a\circ_i\b) = \pm\P(\a)\circ_i \P(\b)
\end{equation}
where the sign is the sign discussed \S\ref{signpar}.
\end{lem}
\begin{proof}
This is simply the observation that first discretizing and then
gluing corresponds to the same combinatorics as first gluing and
then discretizing. This fact  becomes clear if one cuts the graphs
occurring on the right hand side along the simple separating curve
which corresponds to the image of the two glued boundary. If there
are no closed loops and there is not both boundaries are twisted,
the result is immediate. If both the boundaries are twisted this
would yield $0$ in the  open cell cellular chains of $\A$ and hence
on $\arcgraphs$ by definition. We see, however, that we could also
not set them to zero and discretize them. This would again result in
the gluing of discretized twisted arc--graphs, so that the condition
of being twisted in the partitioned case and the non--partitioned
case agree and in both cases the relevant contributions are set to
zero. The same reasoning applies when we erase closed loops. See \S
2.3 of \cite{hoch1} for the definition of the open cell gluing.
\end{proof}

\subsubsection{Self--gluing, Modular Operads and PROPs}
In the above procedures for gluing partitioned arc graphs, we do
not have to assume that the two boundaries we glue actually lie on
the same surface.

\begin{prop}
Allowing self--gluing, the gluing operations of
\S\ref{partgluingpar} turn $\P\arcgraphs(n)$ into a modular
operad. Here the additional modular grading variable ``$g$''
is given by the genus $g$.\footnote{We do not enforce the stability
condition $3g-3+n\geq0$.}

\end{prop}

\subsubsection{Partitioned graphs with $\In$ and $\Out$ markings}
Just like for arc graphs, we can look at partitioned arc graphs
together with a $\Z/2\Z$ marking of their boundaries; viz.\ a
partitioning of their cycles into $\In$ and $\Out$.

\begin{nota}

\label{shortone} To avoid introducing yet other symbols for the
classes of graphs indexing the cells of the different sub--spaces of
$\A$, we simply denote the partitioned graphs by using the prefix
$\P$, e.g.\ $\P\Diioarc$, $\P\iooarc$ and $\P\Arcno$.
\end{nota}

\begin{prop} Restricting the modular operad structure and iterating it by
gluing all $\In$ boundaries to all $\Out$ boundaries of two
collections of elements of $\P\DiA$ imbue a PROP structure on
$\P\DiA$.

Also, similarly, there is a PROP structure ---which we call the
algebraic  PROP structure-- imbued on $\P\DiA$ which is obtained
from the definition of the gluing \S\ref{partgluingpar} by omitting
the final step (1).\footnote{The situation of Step (2) actually
never occurs}.

Finally, the map $\P$ is operadic or in the PROP case is PROPic.
\end{prop}
\begin{proof}
The conditions of associativity are again straightforward if one
cuts the glued surface in two different ways. The operadic
properties are verified as above.
\end{proof}

\subsubsection{Partitioning Angles, grading and the Preservation of the Filtration}
 An angle is a {\em partitioning angle}, if its two
{\em edges}
---that is $e_1=\{f,\imath(f)\}$ and
$e_2=\{\Cyc(f),\imath(\Cyc(f)\}$--- are parallel. In the opposite
case we call it non--partitioning.

With this definition, we can rewrite the grading of $\Parcgraphs$
as given
 graded by  half the number of non--partitioning angles minus one. The number of
non--partitioning angles is precisely the number of edges of the
underlying arc--graph.

Furthermore, it is clear that the degree in the composition in
$\Pcarcgraphs$ goes down by one each time two non--partitioning
angles (other that the outside angles) are glued, as this will
decrease the number of non--partitioning angles by two. In all the
other cases ---gluing partitioning to partitioning and non--partitioning
to partitioning--- the  number of partitioning angles is preserved.
Lastly, erasing closed loops also only decreases the number
of partitioning angles.

Therefore:

\begin{lem} The filtration given by the degree $\leq k$ graded
components on $\Pcarcgraphs$ is respected by the gluings $\circ_i$.
Moreover, the filtration is respected for all gluings in
$\Parcgraphs$ which do not glue an empty boundary to an empty
boundary.
\end{lem}
\qed

\subsubsection{Graded Version}
We also have the same type of statement as in Lemma \ref{Plem} in
the graded case. Recall that we have a grading on $\P\carcgraphs$
and one on $\carcgraphs$, the latter is graded by the dimension of
the cells which is the number of edges -1. This grading is of course
respected by $\P$. Also, in both cases the induced filtration is
respected, so we get an operad structures on the associated graded.

\begin{cor}
The map $\P$ of \ref{Plem} induces an operadic morphism of the
associated graded objects $\P:\Gr\OC(\Arc)\to \Gr\P\carcgraphs$.
Moreover, the same holds true for all sub--operads di-operads or
PROPs whose compositions do not include a glueing of an empty
boundary to an empty boundary.
\end{cor}
\qed

\begin{rmk}
Notice that in this graded version, all the contributions from the
gluing, which involve deleting closed loops are set to zero. This is
true for both sides as deleting a closed loop decreases the grading.
Moreover in  if the condition of step (1) is met, that is if both
boundaries are twisted, the gluing procedure of the algebraic gluing
also decreases the grading, so that the associated graded of the
topological gluing and the algebraic gluing agree.
\end{rmk}

\subsubsection{Angle Marked Partitioned Graphs} We will also
consider the  constructions of the last
paragraphs in the case of angle marked graphs.
  An angle marked partitioned arc graph is a partitioned arc
graph with an angle marking.

We let $\Paarcgraphs(n)^k$ be the angle marked partitioned arc
graphs on a surface with $n+1$  boundaries labelled by $\set{0}{n}$
whose underlying graph has $k+1$ edges and let
$\Paarcgraphs(n)=\bigoplus_k \Paarcgraphs(n)^k$. Then $\Snn$ acts by
permutation on the labels, and we call the collection of $\Snn$
modules $\Paarcgraphs$ simply $\Paarcgraphs$. Again, we use the same
notation for the set and the free Abelian group generated by it.
Also keeping the standard notational conventions, we call the
subset/sub-group of exhaustive partitioned arc graphs
$\Pacarcgraphs$. Also these spaces have a grading by $k$ and hence
an induced filtration by elements of degree $\leq k$.

Analogously to $\P$ there is a partitioning operator for angle
marked arc graphs.
Keeping the notation  $a^p$ for a particularly partitioned
graph as in \S\ref{discpar} set:

\begin{equation}\P(\a,\amark)=\sum
\pm (a^p,\amark_{\a^p})
\end{equation}
 where
$(\amark)^p$ is the angle marking which marks every new partitioning
angle by $1$ and keeps the other markings of $\a$.

 If
$\alpha \in \OC(\DiA)$, where we identify the cell with the
arc graph, then $\a$ has a standard angle marking \cite{hoch1} defined
by marking all outside angles and all angles at the boundaries
 $\In$ by one and the rest, that is the inner angles at the boundaries
$\Out$, by $0$. Likewise there is a standard marking for arc
graphs $\a\in \arcgraphs=CC_*(\A)$
 which was simply the constant marking
by $1$.
For an $\a$ of one of
the two types above, we denote  $\a$ together with its standard
angle marking by
$\alpha^{\angle}\in \aarcgraphs$,
and define

\begin{equation}
\P^{\angle}(\a):=\P(\a^{\angle})
\end{equation}
 Using
the rationale of \cite{hoch1}, we identify a relative cell $\alpha
\in \OCDiarc$ with the angle marked arc graph $\a^{\angle}$ that
labels a cell of $\Ana$, the CW-complex of angle marked arc graphs.

\begin{nota}
Extending the Notation \ref{shortone} we denote the embedding a
class of partitioned graphs into $\Paarcgraphs$ by the prefix $\PA$,
viz.\ $\PA\Diioarc$, $\PA\iooarc$ and $\PA\Arcno$.
\end{nota}

For the gluing we will need a new matching condition.
Given an angle marking, we partition the set of flags $F=F(v)$
at a given boundary into subsets $F=F_1\amalg \dots \amalg F_{k+1}$
where $k$ is the number of markings by $1$ by collecting
together all the flags between which the angle marking is zero.
Notice that $F$ has a linear order and we also think of
the subsets as linearly ordered.

\begin{df}
For {\em two  angle marked partitioned arc graphs} $(\a^p,\amark)$
and $(\a^{\prime p},{\amark}')$ are {\em angle  matched} at the
boundaries $i$ and $i'$ if the number of angles with an angle
marking $1$ agree for these two boundaries. We say that an angle
marked partitioned arc graph is {\em twisted at the boundary $i$} if
the underlying arc graph has this property {\em and} the  each of
the two edges forming the outer angle in the underlying graph has at
least one more parallel edge.

In case two angle marked arc graphs are angle matched at the
boundaries  $i$ and $i'$ there are equally many sets, say $k$, in
the partitions of $F(v_i)$ and $F(v_{i'})$. We say that the graphs
are {\em perfectly} angle matched if $|F(v_i)|>1$ implies
$|F(v_{i'})_{k-i}|=1$ and $|F(v_{i'})|>1$ implies
$|F(v_i)_{k-i}|=1$.
\end{df}

We will give the rigorous combinatorial definition of the gluing
below. Geometrically, we move the edges on the boundaries slightly
apart, if the angle marking between them is $1$. Now the number of
endpoints of these edges on the boundaries which are to be glued
will coincide precisely if they are angle matched. In this case,
we want to identify these vertices. The condition of perfect
matching  ensures that at any given vertex there is at most one
side which has more than one vertex.

\begin{df}
\label{partgluedef} For angle marked arc families $(\a^p,\amark)$
and $(\a^{p \prime},{\amark}')$ which are perfectly angle matched
at the boundaries $i$ and $0$, we define $(\a^p,\amark)\circ_i
(\a^{\prime p},{\amark}')$  as follows. Let $F_{v_0}(\a^p)$ and
$F_{v_i}(\a^{p\prime})$ be the sets of flags at these boundaries.
Add new vertices to all the flags and identify two vertices if two
flags make up an angle with angle marking $0$. Now each of the
sets of vertices obtained form $F_{v_0}(\a^p)$ and
$F_{v_i}(\a^{p\prime})$ comes in a linear order by which we
enumerate them, where we actually enumerate the new vertices
obtained from $F_{v_0}(\a^p)$ in the inverse order, that is last
to first. Now we identify all the vertices with the same number
from the enumeration. We call the flags of  $F_{v_0}(\a^p)$ and
$F_{v_i}(\a^{p\prime})$ flags from different sides, since they lie
on opposite sides of the separating curve that is the image of the
two glued boundaries. Notice that since we are in the perfectly
matched case, there are equally many vertices and if such a vertex
has more than two flags, only one side has more than two flags.
If there are only two flags $\{f_1,f_2\}$ at a vertex, we delete
the vertex and glue the edges, by deleting the two flags
$\{f_1,f_2\}$   and setting
$\imath(\imath_{\a^p}(f_1)):=\imath_{\a^{p\prime}}(f_2)$.  In the
case that there are more flags, say $f$ on one side and
$(f_1,\dots f_l)$ on the other side enumerated in their linear
order, then we duplicate the flag $f$
 $(l-1)$--times and glue the $l$ copies to the $f_i$ in the
 obvious linear order. In this way, we obtain new angles,
namely the angles between the various copies of $f$. We mark all
these angles by $0$. We furthermore forget all the angle markings
at $v_i$ and $v_0$ and retain all other markings.

Again, there are two more steps:

\begin{enumerate}
\item if the arc graphs are twisted at the respective boundaries,
we set the composition  to zero.

\item  if the gluing results  in closed loops, that is embedded
arcs that do not touch any boundary, as edges, set the
contribution to zero.
\end{enumerate}

As above we call this gluing the topological gluing and call the
gluing, which omits step (1) the algebraic gluing.

 {\sc Self--gluing.} The conditions of two boundaries being
perfectly  matched translate in a straightforward fashion to the
case of two boundaries of the same partitioned arc graphs. We
define the self--gluing by the same procedure.
\end{df}

\begin{figure}
\epsfxsize = \textwidth
\epsfbox{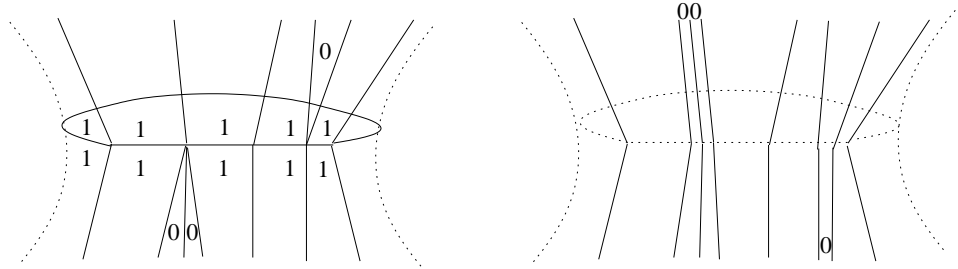}
\caption{\label{gluingex}
An example of a gluing in the perfectly matched case}
\end{figure}

An example of such a gluing is given in Figure \ref{gluingex}.

\begin{rmk}
\label{naturalrmk} The second condition is needed in order to stay
inside the current framework. It is  interesting to remark what
these closed loops correspond to in a different settings. In the
geometric setup they can be understood in terms of general
foliations, see e.g.\ \cite{KP}. In the Hochschild setting, see
\S\ref{hochactionpar} each closed loops contributes a factor that
equals the dimension of the algebra as the trace over the Casimir
element.

The first and second conditions are natural from the topological
point of view, and are necessary if we wish to have an operadic
map from $\carcgraphs$ to $\P\carcgraphs$. Also, see
\S\ref{conclusion} for more comments.
\end{rmk}

\begin{prop}
\label{anglegluingopprop} The gluings defined above turn
$\Paarcgraphs$ bi--graded by the number of boundaries minus one
(as the operadic degree) and the genus (as the genus degree) into
an partial modular operad.

Moreover, this partial  modular operad structure augmented by
setting to zero any thus far not defined gluing is an operad
structure when restricted to $\PA\carcgraphs$ and a di--operad
structure on $\PA\Diioarc$ if one restricts the gluings to gluing
only inputs to outputs. Finally, using consecutive self--gluings
on $\PA\Diioarc$ to glue all ``ins'' to all ``outs'' of
collections of arc graphs,  the partial modular operad structure
induces a PROP structure.

Also, similarly, there is a PROP structure ---which we call the
algebraic PROP structure--- imbued on $\OC\Diioarc$ which is
obtained from the Definition \ref{partgluedef} by omitting the
final step (1).\footnote{Notice that the situation of step (2)
actually never occurs}.
\end{prop}

\begin{proof}
The equivariance with respect to the symmetric group actions is
immediate. The associativity of the partial operations is also
straightforward. Now in both the special cases the condition of
perfect matching is built in and does not change under gluing. In
the first case perfect matching reduces to matching, since all the
angles are marked by $1$. In the case of $\PA\Diioarc$ all the
boundaries which are to be glued are also perfectly matched if they
are angle matched since the $\In$ boundaries are again all marked by
$1$. This condition does not change under gluing, so the result
follows by a straightforward calculation.
\end{proof}

\begin{rmk}
We only defined the gluing in the non--degenerate case. The
general case can be treated in several ways. One is to use the
shuffle combinatorics, like in the definition of gluing for $\Arc$
\cite{KLP,hoch1}. Here one can either average over the occurring
combinatorial types --- or not. Another possibility is to go
outside the present framework and allow arcs running to punctures,
in which case, one simply identifies the vertices and if their
valence is more than two one leaves them in the surface as new
punctures. This is reminiscent of the procedure for open gluing in
\cite{KP}. We shall not need  these considerations
 in the following, but an extension of the gluing
 is of general interest
and deserves further study.
\end{rmk}

\begin{lem}
Let $\a,\b\in \Diioarc$ then
\begin{equation}
\label{anglecompat} \P^{\angle}(\a)\circ_i
\P^{\angle}(\b)=\pm\P^{\angle}(\a\circ_i \b)
\end{equation}
and this map is $\Sn$ equivariant. Thus $\P^{\angle}$ is an
operadic map, if  we use the operadic composition in
$\Gr\OC(\Diioarc)\cong CC_*(\Diioarci)$ on the right hand side.
\end{lem}

\begin{figure}
\epsfxsize = \textwidth
\epsfbox{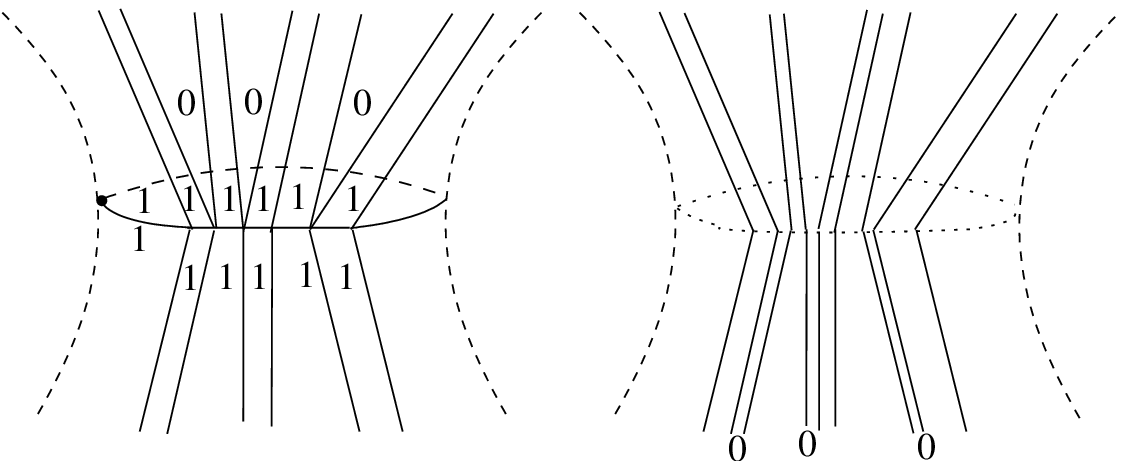} \caption{\label{gluingcompat} An
example of the gluing in $\PA\Diioarc$}
\end{figure}

\begin{proof}
Since in the standard marking, all the angles on the inputs are labelled by one, while the outputs
are labelled by zero, we have the following cases. Two partitioning angles are glued. This does not change the number of non-partitioning angles. An input non-partitioning angle is glued to an output
partitioning angle. In this case the new edges form an non-partitioning angle. Lastly an output non-partitioning angle is glued to a doubled incoming edge which again results in a non-partitioning angle.
So we see that the number of non--partitioning angles is additive.
 (see Figure \ref{gluingcompat} for an example). This means that
 only the graphs of top--degree appear.
 On the other hand it is clear that all possible
 partitioned graphs of $\a\circ\b$ with the maximal number of non-partitioning angles appear.
That is, we obtain exactly the graphs of $\P^{\angle}(\a\circ_i
\b)$.
\end{proof}

\subsection{Graded version for angle marked graphs}

\begin{prop}
$\AP$ is an operadic map between $\Gr\Arcno$
---under the embedding given by marking all angles by one--- and its
 image $\Gr\PA\Arcno$ in $\Gr\Parcgraphs$.
In other words let $\a,\b\in \Gr\Arcno$ then equation
(\ref{anglecompat}) also holds true if we define the left hand
operation to lie
 inside $\Gr\Parcgraphs$.
\end{prop}

\begin{proof}
Notice that for this embedding of the angle marked graphs there
are no restrictions for the partial gluings. In both cases the
boundaries are perfectly angle matched as soon as they are angle
matched.
 In the graded partitioned gluing, we retain the
summands that correspond to the highest possible number of
non-partitioning angles. These in turn correspond to the angles of
the non-partitioned graphs, and thus the conditions imposed on the
gluing correspond to one and another under the map $\PA$.
\end{proof}

\section{The action on a tensor algebra}

\subsection{Actions on the tensor algebra of a vector space}
\label{vectactsec}
As a simpler example than the Hochschild co-chains, we will
consider the action of the various algebraic structures on a
tensor algebra as a warm up in the tradition of
\cite{G}. This type of action is interesting in its own
right. Moreover it is related to the actions on the simplicial
co-chain complexes regarded in \cite{MScosimp}.

\begin{ex}
 Let $V$ be
a vector space and let
$TV=\bigoplus_{n\geq 0}
V^{\otimes n}$ be its tensor algebra.
Let $\TV=\bigoplus_{n\geq 1}
V^{\otimes n}$ be the reduced tensor
algebra. We will now define an action of $\ioarc$ on the tensor algebra.

For each element $\a\in \OC(\ioarc(k,l))$ we consider
$p=(n_1,\dots, n_k)\in P(n,k)$ and consider the summand $\a^p$
where the notation is as in \S \ref{discpar}.
 We let $\a^p$
act on $Hom(\TV^{\otimes k},\TV^{\otimes l})$ as follows. The
operation is defined to be non--zero for $a\in \TV^{\otimes k} $
only if $a\in \bigotimes_{i=1}^k V^{\otimes n_i} \subset
\TV^{\otimes k}$ and zero else. Assume that $a= \bigotimes_{i=1}^k
(a_{i1}\otimes \dots\otimes \a_{in_i})$, we define $\a^p(a)$ as
follows: first decorate the $n$ $In$--boundaries by the given
element of $\TV^{\otimes n}$ using the enumeration of the
boundaries to associate one tensor power of $\TV$ to each boundary
and then associate elements $a_f\in A$
 to each flag $f$
of the arc graph of $\a^p$ whose vertex lies
 on the boundary $j$
by setting $a_f:=a_{i,|f|}$ where $|f|$ is the position of the
flag $f$ in the linear order at the boundary $i$. Now the map
$\imath$ on the level of graphs allows us to associate to $a$ an
element in $\TV^{\otimes m}$ by reading off the elements from the
$\Out$ boundaries as follows: First for the ``out'' boundary $j$
set $\a^p(a)_j:=a_{\imath(f_{j1})}\otimes \dots \otimes
a_{\imath(f_{jm_j})}$  where $f_{j1}, \dots, f_{jm_j}$ are the
flags at the boundary $j$ in their linear order. Then set $\a^p(
a)=\bigotimes_j \a^p( a)_j\in \TV^{\otimes m}$.

It is clear that the action $\a(a):=\P(\a)(a)$ is operadic. Moreover
embedding the open cell complex of $\OC(\ioarc)$ into
$\OC(\Diioarc)$ and using the algebraic PROP structure on $\Diioarc$
this action becomes a PROP action.
\end{ex}

\begin{prop}
\label{vectaction} $\TV$ is an algebra over the PROP
$\OC(\ioarc)\subset \OC(\Diioarc) $ endowed with the algebraic
PROP structure. Moreover this action extends to an action of
$\OC(\Diioarc)$ on $TV$.
\end{prop}
\begin{proof}
The first statement is straightforward. In the non-exhaustive
case, we decorate the empty boundary components by elements of
$k=V^0$ that is after using the imposed multilinear properties we
simply insert $1\in k$ in the above calculations. The effect is
that the empty boundaries are simply ignored using the isomorphism
$V\otimes_k k \cong V$.
\end{proof}

\begin{rmk}
Notice that this action involves only the boundary data and thus only the ribbon graph and not its specific embedding. Thus the action actually factors though the map $\Loop$ of \cite{KLP} and will give an
example of an action of the stabilized Arc operad \cite{Ribbon}.
\end{rmk}

\subsection{The graded version}

\subsubsection{The sub--operad $\Shuff$ of $\CHom(TV)$.}
Let $V$ be a  vector space and let $\CHom(TV)(n):=Hom(TV^{\otimes
n},TV)$ be the endomorphism operad of $TV$. It is well known that
$TV$ is an algebra for the multiplication $\mu_{\otimes}:=\otimes$
and a  co-algebra for the co-multiplication $\Delta:TV\to TV\otimes
TV$ given by $\aoa{1}{n}\mapsto \sum_i (\aoa{1}{i})\otimes
(\aoa{i+1}{n})$. Also, we have the $\Sn$ action on $TV^{\otimes n}$
permuting the factors of $TV$. These three basic sets of operations
generate a suboperad of $\CHom(TV)$ which we would like to call
$\Shuff$. It can be thought of as something like the brace algebra
which is a natural subalgebra of operations on the Hochschild
co--chains on an associative algebra given by the natural operations
in that setting. In this spirit $\Shuff$
 are the natural operations on a tensor algebra.

\begin{nota}
\label{deltanotation}
We will use the following notation:
$\Delta^{l}:TV\to TV^{\otimes l+1}$ is the iteration of $\Delta$ given by
$(\Delta \otimes id^{\otimes l} )\circ (\Delta \otimes id^{\otimes l-1})
\circ \dots \circ (\Delta \otimes id) \circ \Delta$.
\end{nota}

We wish to point out that any element in $\Shuff$
can be uniquely written as a sum of elements of the type
$\mu_{\otimes}^{\otimes  \sum n_i +(n-1)}
 \circ \sigma\circ \bigotimes_i \Delta^{n_i}$ where one
first uses the coproduct $n_i$ times
on the $i$-th factor of $TV$ in the product $TV^{\otimes n}$,
 then one uses a permutation $\sigma$ on the resulting factors
of $TV$  and finally
one multiplies them all together.

$\Shuff$ is graded by $\sum_i (n_i+1)-1$.
The composition does not respect this grading, but it does respect
the induced filtration, and hence we get an induced operad
structure on the associated graded $\Gr\Shuff$. The same type
considerations {\em mutatis mutandis}
 hold true for the PROP
$Hom(TV^{\otimes n},TV^{\otimes m})$ and we use $\Shuff(n,m)$
to be the appropriate sub--PROP.

The PROP action by $\OC(\Diioarc)(n,m)$ is then contained in
$\Shuff(n,m)$ and  it is actually easy to see that the image is
precisely $\Shuff(n,m)$ by using the normal form above. The
grading then  corresponds to the  number of angles on the inputs
minus one which is the number of edges minus one that is the
dimension of the cell considered in $\Diioarci$, so we see that
these two gradings are compatible.

\begin{prop}
The PROP action of $\OC(\Diioarci)$ on $TV$ has its image in
$\Shuff$ and passes to the associated graded $\Gr\Shuff$, so that
there is an operadic (PROPic) map from $\Gr\OC(\Diioarc)\simeq
CC_*(\Diioarci)\rightarrow \Gr\Shuff$.
\end{prop}\qed

\begin{rmk}
The condition of genus zero also has a nice algebraic
counterpart in this setting and that is the condition
 that the permutation is only a shuffle.
\end{rmk}

\begin{ex}
An example where the grading is not respected
occurs when one considers $\Delta \circ \mu_{\otimes}: TV\otimes TV\to TV
\otimes TV$.
The generic number of shuffles will be $1$, but there
will be a summand corresponding to $id\otimes id$ which
will require no shuffle.

This example is very instructive, since it is this sort of behavior
which is not Frobenius that is very characteristic for our actions
and their associated graded ones.

We can decompose
 \begin{equation}
\Delta \circ \mu_{\otimes}=
(id \otimes \mu_{\otimes})\circ (\Delta \otimes id) + id\otimes id
+ (\mu_{\otimes}\otimes id)\circ (id\otimes \Delta)
\end{equation}
As well known this means that $TV$ as an algebra and co--algebra
is not Frobenius but rather has the operations usually forming the other
side of the Frobenius equation
as summands of the operation $\Delta \circ \mu_{\otimes}$;
plus there is one more summand of lower degree, namely $id\otimes id$.
The associated graded will project out this term.

\end{ex}

In order to go  beyond the PROPic setting of  $\ioarc$
and $\Diioarc$ in the setting of the example of Proposition \ref{vectaction},
we will need a pairing. The exact axiomatic setup for this is
given in the next section.

\section{Operadic Correlation Functions}

In this section, we introduce operadic correlation functions,
which can be thought of as the generalization of an algebra over a
cyclic operad to the $dg$--setting. In order to get to the main
definition, we first set up some notation.

Given a pair $(A,C)$ where
$A$ is a vector space and
$C=\sum c^{(1)}\otimes c^{(2)}\in A\otimes A$
we define the following operations
\begin{equation}
\label{cmult}
\circ_i:Hom(A^{\otimes n+1},k)\otimes Hom(A^{\otimes m+1},k)
\to Hom(A^{\otimes n+m},k)\\
\end{equation}
where for $\phi \in Hom(A^{\otimes n+1},k)$ and
$\psi \in Hom(A^{\otimes m+1},k)$
\begin{multline}
\phi\circ_i \psi(\aoa{1}{n+m})=\\\sum \phi(\aoa{1}{i-1}
\otimes c^{(1)} \otimes \aoa{i+m}{m+n})\psi(c^{(2)}\otimes
\aoa{i}{i+m-1})
\end{multline}

\begin{df} A set of operadic correlation function for a cyclic linear
operad $\O$ is a tuple $(A,C,\{Y_{n}\})$ where $A$ is a vector
space,  $C=\sum c^{(1)}\otimes c^{(2)}\in A\otimes A$ is a fixed
element and  $Y_{n+1}:\O(n)\to Hom(A^{\otimes n+1},k)$ is a set of
multi--linear maps.  The maps $\{Y_{n}\}$ should be $\Snn$
equivariant and for $op_n\in \O(n), op_m \in \O(m)$

\begin{equation}
\label{corrcompat}
Y_{n+m}(op_n\circ_i op_m)=
Y_{n+1}(op_n) \circ_iY_{m+1} (op_m)
\end{equation}
where the $\circ_i$ on the left is the multiplication of equation
(\ref{cmult}) for the pair $(A,C)$.

We call the data $(A,\{Y_{n}\})$ of an algebra and the $\Snn$
equivariant maps correlation functions or simply correlators for
$\O$.
\end{df}

\subsection{Correlators for algebras over cyclic operads}
An example is given by an algebra over a cyclic operad.
Recall that this a triple $(A,\pair, \{\rho_n\})$
where $A$ is a vector space, $\pair$ is a non--degenerate
bi--linear pairing and $\rho_{n}:\O(n)\to Hom(A^{\otimes n},A)$
are multilinear maps also called correlators that satisfy
\begin{itemize}
\item[i)] $\rho(op_n\circ_i op_m)=\rho(op_n)\circ_i\rho(op_m)$
where $\circ_i$ is the substitution in the i-th variable.
\item[ii)] The induced maps $Y_{n+1}:\O(n)\to Hom(A^{\otimes n+1},k)$
given by
\begin{equation}
\label{opcorr}
Y_{n+1}(op_n)(\aoa{0}{n}):=\la a_0,\rho(op_n)(\aoa{1}{n}) \ra
\end{equation}
are $\Snn$ equivariant.
\end{itemize}

\begin{nota}
Given a finite dimensional vector space  $A$ with a non-degenerate
pairing $\pair=\eta\in \check A\otimes \check A$, let $C\in
A\otimes A$ be dual to $\eta$ under the isomorphism induced by the
pairing and call it the Casimir element. It has the following
explicit expression: Let $e_i$ be a basis of $V$, let
$\eta_{ij}:=\la e_i,e_j\ra$ be the matrix of the metric and let
$\eta^{ij}$ be the inverse matrix. Then
$C=\sum_{ij}e_i\eta^{ij}\otimes e_j$.
\end{nota}

\begin{lem}
The assignment $(A,\pair,\{\rho_n\}) \mapsto (A,C,\{Y_n\})$ where
the $Y_n$ are defined as in equation (\ref{opcorr}) gives a  1--1
correspondence between the algebra structure over a cyclic operad
and
 operadic correlation functions, which is
functorial.
\end{lem}

\begin{proof}
We have defined the map in one direction. To give the inverse
map, we set
\begin{equation}
\rho_n(\aoa{1}{n}):=\sum Y_{n+1}(c^{(1)},\aoa{1}{n})c^{(2)}
\end{equation}
A direct calculation verifies that these assignments are inverse
to each other. The compatibility of the $\Snn$ operations is
manifest. Finally, it is clear that this construction is functorial
for maps of cyclic operads and maps of algebras with a non--degenerate
pairing.
\end{proof}

\begin{ex}
\label{tvshuffex}
We can now generalize the example of Proposition \ref{vectaction}.
For this fix a non-degenerate
symmetric pairing $\pair$ on $V$.
Then
this pairing induces a symmetric non-degenerate pairing on $TV$, so
that we can consider $TV$ (and $\TV$) as a candidate of an algebra
over a cyclic operad. We give this  structure via operadic
correlations functions in analogy to the operations of
\S\ref{vectactsec}.

Let $\a\in (\carcgraphs)^k(n)$ and let $p\in P(m,k)$. We define
an $Y(\a^p):\TV^{\otimes n}\to k$ as follows. Let $n_i:=\#\{$flags
at the boundary $i\}$ of the arc graph of $\a^p$.

If  $ a\notin \bigotimes_{i=1}^k V^{\otimes n_i} \subset
\TV^{\otimes n}$, we set $\a^p( a)=0$. And if $ a=\bigotimes_{i=1}^k
(a_{i1}\otimes \dots \otimes \a_{in_i})$ then we associate to each
of the flags $f$ of the arc graph of $\a^p$ incident to the boundary
$i$ the element $a_f:=a_{i,|f|}$, where again $|f|$ is the position
of the vertex in its linear order. Let $E$ be the set of edges of
the arc graph of $\a^p$  then we define $Y(\a^p)(
a):=\prod_{e=\{f,\imath(f)\}\in E}\la f,\imath(f)\ra$. Now it is
again straightforward to check that this defines operadic
correlation functions for $\P\carcgraphs$ with the algebraic gluings
and in the case that the final step (2) of \S\ref{partgluingpar} is
not applicable. One can then extend to operations of $\Parcgraphs$
and hence $\arcgraphs$ on $TV$ by again using the isomorphisms
$V\otimes_k k\simeq V$ to ``decorate'' the empty boundaries with
copies of $k$.
\end{ex}

\subsection{Operadic correlations functions with values in a
twisted $\CHom$ operad}

\begin{df}
Let $(A,\pair,\{Y_{n}\})$ be as above. And let $\H=\{\H(n)\}$ with
$\H(n)\subset \Hom(A^{\otimes n},A)$ as $k$--modules be an operad
where the $\Sn$ action is the usual action, but the operad structure
is {\em not necessarily} the induced operad structure. Furthermore
assume that $\rho_{Y_{n+1}}\in \H(n)$. We say that the $\{Y_{n}\})$
are operadic correlation functions for $\O$ with values in $\H$ if
the maps $\rho$ are operadic maps from $\O$ to $\H$. We will also
say that we get an action of $\O$ with values in $\H$.
\end{df}

\subsubsection{Signs} As in the case of the Deligne conjecture one
twist which we have to use is dictated by picking sign rules. In the
case of Deligne's conjecture this could be done by mapping to the
brace operad $\Brace$ (see e.g.\ \cite{del}) or by twisting the
operad $\CHom$ by lines of degree $1$ as in \S\ref{signpar} (see
e.g.\ \cite{KS}). In what follows, our actions will take values on
operads that are naturally graded {\em and moreover} we will
identify the grading with the geometric grading by  e.g.\ the number
of edges or the number of angles etc.. The signs will then
automatically match up, if we use the procedure of \S\ref{signpar}
at the same time for {\em both} the graph side and the $\CHom$ side,
i.e.\ for the operad $\H$. In fact, this approach unifies the two
sign conventions mentioned above on the subspace of operations
corresponding to $\Lintree_{cp}$.

\subsection{Actions on the tensor algebra of a vector space}
Assume now that $(V,\pair)$ is a finite dimensional vector space
with a non--degenerate pairing.

\begin{df}
We let $\Modshuff(n)\subset \Hom(TV)(n)$ be the image of all the operations of
$\P\arcgraphs(n)$, by considering $0$ as ``out''.
\end{df}

\begin{nota}
\label{etanota} We write $\eta^m$ for the map $V^{\otimes 2m}\to
k$, $\eta(\aoa{1}{2m})=\prod_{i=1}^m \eta(a_i,a_{2m-i})$.
\end{nota}

\begin{prop}
\label{modshuffprop}
After dualizing
to obtain elements in $Hom(TV^{\otimes n+1},k)$ any element
in $\Modshuff(n)$
can be written uniquely as

\begin{equation}
\bigotimes_{j=1}^{\frac{1}{2}\sum_i (n_i+1)}\eta^{m_j}\circ \sigma \circ
\bigotimes_{i=0}^n \Delta^{n_i}
\end{equation}
where  $\sigma$ is a permutation of the $\sum_i (n_i+1)$--factors
of $TV$,
and we used the Notation \ref{deltanotation}.
Set $l=\frac{1}{2}\sum_i (n_i+1)-1$ then $\Modshuff$ is graded by
$l$. Moreover the composition in $\Modshuff$ respects
the induced filtration of elements of degree $\leq l$. Lastly,
 the decomposition identifies $\Modshuff$
with the subspace of $\Hom(TV)$ obtained by dualization for the
subspace generated by the coproduct, permutations and $\eta$ in
$\bigoplus_n Hom(TV^{\otimes n+1},k)$. That is we obtain
correlation functions with values in $\Modshuff$.
\end{prop}

\begin{proof}
The first statement is clear by the definition of $\Modshuff$ as the image.
The last statement is also straightforward, by arranging the
operation in the specified order. On the other hand it is easy to give
the arc graph in $\P\arcgraphs(n)$ by drawing one arc for each factor
of $\eta$ with the incidence relations given by $\sigma$.
This identifies the two subspaces. In this identification there is one
factor of $\Delta$ for each non-partitioning inner angle.
The last claim, that the operations respect the filtrations
is clear after identifying $k$ with the dimension of the cell, that is
the number of edges minus one, of the underlying graph for the operation.
The mentioned equality follows from the combinatorial
identity $|\angle_{\rm inner}|+|\angle{\rm outer}|=|{\rm Flags}|
=2 |{\rm edges}|$.
\end{proof}

\begin{prop}
\label{tvshuffprop}
 For any vector space $V$ with a non-degenerate pairing,
$TV$ is an algebra over the algebraic PROP $\P\iooarc$. Moreover
this action passes to the associated graded and gives an action of
$CC_*(\A)$ with values in $\Gr\Modshuff$.
\end{prop}

\begin{proof}
We extend the definition of correlators above to correlators for
$A$ on $TV$ as in \S\ref{vectaction} above. We use the Casimir
element to dualize and thus we only have to show that the
resulting structure is that of an algebra over a cyclic operad.
Again the $\Snn$ equivariance is manifest. After dualizing the
gluing on the flags in the operadic composition turns into the
identity map $id: A\to A$ so that  indeed the gluing $\circ_i$ on
$\P\iooarc$ maps to insertion at the $i$-th place in $\Hom(TV)$.
Dealing with the extra steps (1) and (2), we see that on the side
of $\Hom(TV)$ they would not yield zero. In $\P\iooarc$ closed
loops cannot appear hence step (2) is avoided. Using the algebraic
gluing, we get agreement for the two operations. In the graded
case, the contributions of (1) and (2) are projected away on the
side of $\Modshuff$ as is the case in $CC_*(\A)$ where these
contributions come from lower degree cells, which are again
projected out. Moreover as mentioned above the grading in both
cases is by the number of arcs of the underlying graphs -1.
\end{proof}

\subsection{Correlators for $dg$--algebras}

Let $(V,d)$ be a complex  whose homology algebra $H:=H(V,d)$ has a
non--degenerate pairing $\pair$. Let $C$ be the Casimir element of
$\pair$. Let $Z=\ker(d)$ and let $s$ be a section of the
projection $Z\to H$.

\begin{df}
If $(V,d)$ is a complex, we call a set of  correlation functions for
$V$  operadic chain level correlation functions if they are operadic
correlation functions for $Z(A)=\ker(d)$.
\end{df}

\begin{prop}
\label{dgcorrelatorprop} Let $\O:=\{\O(n)\}$ be a cyclic operad, and
let $(V,d,\pair)$ and $s$ be as above. Let $Z=\ker(d)$ and let
$Y_{n+1}:\O(n)\to Hom(V^{\otimes n+1},k)$ be a collection of $\Snn$
equivariant maps
 whose
values only depend on the classes in $H$, that is if $[a_i]=[b_i]\in
H$ where $a \mapsto [a]$ is the projection map $Z\to H$,  then
$Y_n(\aoa{1}{n})=Y_n(\bob{1}{n})$. Furthermore assume that the
induced maps $\overline Y_{n+1}: \O(n)\to Hom(H^{\otimes n+1},k)$
are operadic correlation functions for the cyclic operad, then for
any section $s$ of the projection map $Z\to H$ the collection
$\{Y_n\}$ is a set of operadic correlation functions for
$(Z,(s\otimes s)(C))$.
\end{prop}

\begin{proof}
Straightforward.
\end{proof}

An example is given by adapting Example \ref{tvshuffex}
to the current setting.

\begin{thm}
 For a complex $(V,d)$ over $k$ with a pairing $\pair$ which
satisfies
\begin{itemize}\item[i)]
 $\forall \, v,w \in V: \la dv,w\ra=-\la w,dv\ra$ and
\item[ii)] the induced pairing on $H=H(V,d)$ is non--degenerate,
\end{itemize}
 the formulas of Example \ref{tvshuffex}extended as in Proposition \ref{tvshuffprop}
 yield
correlation functions on the tensor algebra $TV$ which are operadic
on $TZ$ that is operadic chain level correlation functions.
\end{thm}
\qed

\begin{df}
\label{qfidef} A quasi--Frobenius algebra is a triple
$(A,d,\pair)$ where  $(A,d)$ is a unital $dg$--algebra  whose
homology algebra $H:=H(A,d)$ is finite dimensional and has a
non--degenerate pairing $\pair$ and is a Frobenius algebra for
this pairing. A quasi--Frobenius algebra with an integral is a
triple $(A,d,\int)$ where $\int: A\to k$ is a linear map such that
\begin{itemize}
\item[i)]
$\forall a\in A:\int da=0$
\item [ii)] $(A,d,\pair)$ is a quasi--Frobenius
algebra,where $\la a, b\ra:=\int ab$.
The cocycles of a quasi--Frobenius algebra with an integral are the
subalgebra $Z=\ker(d)\subset A$ of the algebra above.
\end{itemize}
\end{df}

A natural example of a quasi--Frobenius algebra with an integral is
 $A=C^*(M)$, the co-chains of a compact manifold $M$.
\begin{ex}
\label{cyclicex}
Let $\Cyclic:=\Ass[1]$ be the cyclic operad obtained by shifting the associative operad $\Ass$
by $1$ that is $\Cyclic(n)$ is the
 permutation representation $\Snn$ on $\set{0}{n}$
 with $\Sn$ acting on $\set{1}{n}$.

Let $(A,d)$ be an associative $dg$--algebra with an integral,
i.e.\ a function $\int:A\rightarrow k$ which satisfies $\int ab=
(-1)^{\deg(a)\deg(b)}\int ba$ and Set $\la a,b\ra =\int ab$ and
assume $\pair$ is non--degenerate on $H=H(A,d)$.

For $\sigma_n\in \Sn$ viewed as a generator of $\Cyclic(n-1)$ we
define
\begin{equation}
Y_{n+1}(\sigma_n)(\aoa{0}{n-1}):=\pm \int \ata{\sigma(0)}{\sigma(n-1)}
\end{equation}
where $\pm$ is the sign of the permutation of the elements $a_i$.
Then  by the proposition above
this is a set of operadic correlation functions.

Fixing a section $s$ of the map $Z\to H, a \mapsto [a]$ the operations
\begin{equation}
\rho_n(id_{n+1})(\aoa{1}{n})=\pm s([a_1] \dots [a_n])
\end{equation}
yield operadic correlation functions.

In the case that $A=C^*(M)$ we see that we recover the cup product up to homotopy.
This is enough to characterize the cup product of two closed
 co-chains inside any integral.
\end{ex}

\subsection{Polygon correlation functions}
Other examples of operadic correlation functions come from operads
of polygons. It is this type of example which we generalize to
obtain the correlators for  $\P\Anarc$.

Let $p_n$ be the regular $n$-gon and denote its sides by
$sides(p_n)$. We let $\Poly(n)=$ the free Abelian group generated
by $\{\lab:sides(p_n)\rightarrow \{0, \dots, n-1\}$. Similarly,
let $P_{2n}$ be the regular $2n$-gon with a fixed choice of a
preferred set of $n$ non-intersecting sides, which we call
$Sides(P_{2n})$. Then $\Poly_2(n)=$ the free Abelian group
generated by $\{\lab:Sides(P_{2n})\rightarrow\{0, \dots,n-1\}\}$.
The operad structure $\circ_i$ on $\Poly(n)$ is given by gluing
the polygons along the sides marked by $0$ and $i$, respectively,
and deleting the image of the glued side, which is  diagonal in
the glued object. In $\Poly_2(n)$ we also merge the two pairs of
non--labelled sides on the two sides of the deleted diagonal. The
$\Snn$ action is given by permuting the labels on the labelled
sides.

\begin{rmk}
\label{polycyclic}
These operads are different incarnations of the operad $\Cyclic$.
The map from $\Poly$ to trees is just given by marking the center
of the polygon and the middle of the sides by a vertex, connecting
the vertices of the sides to the center vertex and carrying over the labelling.
This yields an isomorphism to the usual pictorial way of defining $\Cyclic$ in
terms of planar corollas (see e.g.\ \cite{woods}).
An isomorphism from $\Poly$ to $\Poly_2$ is given by
blowing up the vertices of $p_n$ to sides and choosing the
original sides to be the preferred set of $n$ sides of $P_{2n}$.
The inverse is given by contracting the non-preferred sides.
\end{rmk}

Let $(A,d,\int)$ be a quasi--Frobenius algebra with an integral.
 We define correlation functions  as

\begin{eqnarray}
\label{polycoreq}
&&Y(\lab:sides(p_n))\rightarrow
\set{0}{n-1})(\aoa{\lab^{-1}(0)}{\lab^{-1}(n-1)})\nn\\&&\qquad=
\int
\ata{\lab^{-1}(0)}{\lab^{-1}(n-1)} \\
&&Y(\lab:Sides(P_{2n}))(\aoa{\lab^{-1}(0)}{\lab^{-1}(n-1)})\nn\\
&&\qquad=
 \int \ata{\lab^{-1}(0)}{\lab^{-1}(n-1)}
\end{eqnarray}
here and in the following we will frequently drop the subscripts on the $Y_n$ since
they can be deduced from the expression.

\begin{lem}
Given a quasi--Frobenius algebra with an integral $(A,d,\int)$,
the equations above define operadic correlation functions for
$\Poly$ and $\Poly_2$.
\end{lem}
\begin{proof}
This follows either directly by Proposition \ref{dgcorrelatorprop}
or  by Remark \ref{polycyclic} and the Example \ref{cyclicex}.
\end{proof}

\subsubsection{$\Ainf$-algebras and polygons with diagonals}
Another example  of operadic correlation functions comes from
$\Ainf$ algebras. Let $A$  be an $\Ainf$--algebra with
multiplications $\mu_n: A^{\otimes n}\rightarrow A$ and set
$d=\mu_1$.

 Define  $\Polydiag$, the operad of
polygons with diagonals, as follows: $\Polydiag(n)$ is the operad
generated by the free $\Snn$ module generated by the set
$\{diag(p_{n+1})\}:=\{$polygonal decompositions of the abstract
planar $n+1$-gon whose sides are cyclically labelled
$\set{0}{n}\}$. $\Snn$ acts by permutations on the labels. We call
this set $diag$ since the decomposition amounts to choosing
several non-intersecting diagonals. The operad structure on
$\Polydiag$ is given by gluing the polygons along the indicated
sides and keeping the image of the glued sides as a diagonal in
the glued polygon. There is a natural $dg$--structure on this
operad whose differential is the sum of elements obtained by
inserting different diagonals with the appropriate sign. The sign
is determined considering $L^{\otimes \{\rm{polygonal \,
regions}\}}$ as described in \S\ref{signpar}. It is clear that as
a collection of $\Snn$-modules $\Polydiag$ is generated by the
elements $o_n$ where $o_n$ is the cyclically labelled abstract
$n+1$-gon with no diagonals.

\begin{rmk}
Again this operad is just a re-writing of an old, familiar operad.
This time it is Stasheff's $\Ainf$ operad in its tree description,
see e.g. \cite{MSS}. The isomorphism is given by considering the
dual tree of the polygonal decomposition. This tree has one vertex
for each polygonal regions and one for each of the labelled sides.
The edges are given by connecting two vertices if they have a common
diagonal or if the labelled side is a side of the polygonal region.
The vertices corresponding to the labelled sides are exactly the
vertices of valence one and are naturally labelled. This graph is
easily seen to be a tree and the cyclic order induced by the cyclic
order on the sides of the polygonal regions induced by the
orientation of the plane makes this tree into a planar tree. Fixing
the root to be the vertex corresponding to the side labelled by
zero, we obtain a planar planted tree whose leaves (viz.\ non--root
vertices of valence one)  are labelled. The image of the
differential under this correspondence will contract edges.
\end{rmk}
 Assume that $A$ is finite dimensional and has an pairing
non--degenerate $\pair$ and let $C$ be the Casimir element of the
pairing, then set

\begin{equation}
\label{snncond} Y(o_n)(\aoa{0}{n})=\la a_0,\mu_n(\ata{1}{n}) \ra
\end{equation}

\begin{lem}
Given an $\Ainf$-algebra $A$ over $k$ with an
 a non--degenerate pairing $\pair$
such that the correlation functions (\ref{snncond}) above are
$\Snn$ equivariant, extend the definition of $Y$ operadically, by
using the equation (\ref{corrcompat}) recursively as a definition.
Then this extension gives  operadic correlation functions for the
cyclic operad $\Polydiag$. These correlation functions are even
compatible with the $\Ainf$ differential.
\end{lem}

\begin{proof}
The fact that we obtain operadic correlation functions is true by
construction. The fact that they are compatible with the
$\Ainf$--differential then follows from the fact that in  standard
tree depiction of the $\Ainf$--algebra (see e.g.\ \cite{MSS})  the
differential adds edges, which dually corresponds to inserting a
diagonal as explained in the remark above.
\end{proof}

\section{Correlators for $\PA\A$ on the Hochschild co--chains
of a Frobenius algebra}

The natural operations on the tensor algebra $TA$ of an algebra
$A$ are the ones generated by the multiplication and
co--multiplication $\mu_{\otimes}$ and $\Delta$ of $TA$ as well as
the permutations and the multiplication $\mu_A:A\otimes A\to A$.
In order to incorporate the later operations into the picture, we
will have to modify the correlators of \ref{tvshuffex} a little
and introduce operations which act ``internally'', that is
operations which are associated to complementary regions or dually
in the case of $\Arcno$ at the vertices of the dual graph.

\subsection{Graph correlation functions aka.\  Feynman rules}

\begin{ex} We will recall how
to define correlation functions for ribbon graphs,
(see e.g.\ \cite{Keuro,KM}), by using so-called Feynman rules.
Although our action will be slightly different, the underlying
principle is similar and this easier example will be instructive.

Let $\G$ be a ribbon graph with vertices of valence at least 3.
Let $A$ be an algebra with a non-degenerate pairing, which gives
an isomorphism of $A$ with its dual $\check A$. Let $C$ be the
Casimir element. Let $\phi:V_{\G} \rightarrow \CHom(A,A)$ be a map
which

\begin{itemize}
\item[i)]
 preserves degree, i.e.\ $\phi(v)\in
\CHom(A,A)(\val(v-1))\simeq \Hom(A^{\otimes \val(v)},k)$.

\item[ii)] has a cyclic image, i.e.\
$$\forall vf\in V_{\G} :
\phi(v)(a_0 \otimes \dots\otimes a_n)= \pm\phi(v)(a_1\otimes \dots
\otimes a_n \otimes a_0)$$
where $\pm$ is the permutation
super-sign.
\end{itemize}

Set
\begin{equation}
\label{YG} Y(\G)(\phi):= (\bigotimes_{v\in V_{\G}(\t)}\phi(v))
(\bigotimes_{e\in E_{\G}} C)\in k
\end{equation}

\end{ex}

\begin{rmk}
Let $\mathcal{G}$ be a class of graphs, e.g.\
 the set of all planar trees or ribbon graphs. Functions
of the type $\phi:\mathcal{G}\rightarrow \CHom(A,A)$ are sometimes
called Feynman rules if the map $\phi$ is expressible in local
data of the graph, that is in terms of the flags at each vertex
and the edges.

Graphs, which have external vertices, that is vertices with only
one adjacent flag,
 are treated in one of the following ways. One can either partition the
external vertices into a set $In$ and a set $Out$ say of
cardinalities $p$ and $q$ and view $Y(\G)^{p,q}\in Hom(A^{\otimes
p},A^{\otimes q})$ or  view all vertices as inputs $Y(\G)\in
Hom(A^{\otimes p+q},k)$. An example is given by associating
elements of $A$ to the flags of the external vertices of a ribbon
graph with tails   and plugging in copies of $C$ into the
operations by decorating the internal edges and then contracting
like in (\ref{YG})
\end{rmk}

\begin{ex}
\label{Kex} An important example for Feynman rules $\phi$ for
ribbon graphs $\G$ is given by using an $A_{\infty}$--algebra
\cite{Keuro}. Here $\phi(v)=\phi(\val(v))=\mu_n$ where now
$\mu_n:A^{\otimes n}\rightarrow A$ is one of the structure maps of
the $A_{\infty}$--algebra.
\end{ex}

\begin{ex}
In \cite{cyclic} we also used Feynman rules for trees with external vertices.
They are related to the operations which we describe below, by moving
to the dual graph and then to the intersection graph.
They are hence decidedly different for the Example \ref{Kex} above, which
directly deals with ribbon graphs.
\end{ex}

\subsection{$\Anarc$ correlation functions}
The idea of how to obtain the correlation functions for the tensor
algebra is very nice in the $\Arc$ picture where it is based on
the polygon picture. This polygon picture can be thought of as an
IRF (interaction 'round a face) picture for a grid on a surface
which is dual to the ribbon picture. For this we would
modify the arc graph by moving the arcs a little bit apart as described in \S\ref{drawpar}.
Then the complementary regions of partitioned
 quasi--filling arc--graphs
$\P\qfarcgraphs$ are $2k$-gons whose sides alternatingly correspond
to arcs and pieces of the boundary. The pieces of the boundary
correspond to the angles of the graph and of course any polygonal
region corresponds to a cycle of the arc graph. If the graph
$\alpha^p$ has an angle marking, then the sides of the polygons
corresponding to the boundaries will also be marked. We fix the
following notation. For an angle marked partitioned arc graph
$\a^p\in \PA\qfarcgraphs$ let $Poly(\alpha^p)$ be the set of
polygons given by the complementary regions of $\alpha^p$ when
treated as above.
 For $\pi\in
Poly(\alpha^p)$, let $Sides'$ be the sides
corresponding to the angles which are marked by $1$ and
$Sides'(\alpha^p)$ be the union of all of these sides.
If we denote $\angle^+(\G)=(\amark)^{-1}(1)$
 there is a natural bijection between $\angle^+(\a^p)$ and
$Sides'(\a^p)$.

\subsubsection{Correlation functions on the tensor algebra of an algebra}
Fix an algebra $A$ with a cyclic trace, i.e.\ a map $\int:
A\rightarrow k$ which satisfies $\int a_1 \dots a_n= \pm \int a_n
a_1\dots a_{n-1}$ where $\pm$ is the standard  sign.

Now for $\pi \in Poly(\a^p)$ set
\begin{equation}
\label{picor}
 Y(\pi)(\bigotimes_{s\in Sides'(\pi)} a_s)=\int \prod_{s\in
Sides'(\pi)} a_s
\end{equation}
Notice that we only have a cyclic order for the sides of the
polygon, but   $\int$ is (super)-invariant under cyclic
permutations, so that if we think of the tensor product and the
product as indexed by sets  (\ref{picor}) it is well defined.

For an angle marked partitioned arc family $\a^p$ set

\begin{equation}
\label{partcorreq}
Y(\alpha^p)(\bigotimes_{s\in (\amark)^{-1}(1)}
a_s)=\bigotimes_{\pi\in Poly(\alpha^p)}Y(\pi)(\bigotimes_{s\in
Sides'(\pi)} a_s)
\end{equation}
where we  used the identification of the set
$Sides'(\a^p)=
\amalg_{\pi \in Poly(\a^p)} Sides'(\pi)$ with $\angle^+(\a^p)$.
Since for each $\a^p\in \PA\carcgraphs(n)$
 the set of all flags has a linear order,
we can think of $Y(\a^p)$ as a map $A^{\otimes
|F(\a^p)|}=\bigotimes_{i=1}^n A^{\otimes |F(v_i)|} \to k$ and
furthermore as a map to $TA^{\otimes n}\to k$ by letting it be
equal to equation (\ref{partcorreq}) as a map from
$\bigotimes_{i=1}^n A^{\otimes |F(v_i)|} \subset TA^{\otimes n}$
and setting it to zero outside of this subspace.

Extending linearly, for an angle marked arc family $\alpha\in \Anarc$,
we finally define
\begin{equation}
\label{ydef}
Y(\alpha):=Y(\P(\a))
\end{equation}

\subsubsection{Correlators for the
Hochschild co-chains of a Frobenius algebra }
\label{hochactionpar}

Let $A$ be an algebra and let $\CH^n(A,A)=Hom(A^{\otimes n}, A)$
be the Hochschild cochain complex of $A$. We  denote the cyclic
cochain complex by $CC^n(A,k)= Hom(A^{\otimes n+1},k)$. Then one
has a canonical isomorphism of $\CC^*(A)\cong \CH^*(A,\check A)$
as complexes and hence also $\HC^*(A)\cong \HH^*(A,A)$ where $\HC$
is Connes' cyclic cohomology and $\HH$ is the Hochschild
cohomology.

\begin{lem}
For any
Frobenius algebra $(A,\pair)$, we have canonical isomorphisms
$\CC^*(A)\cong \CH^*(A,\check A)\cong \CH^*(A,A)$
and $\HC^*(A)\cong \HH^*(A,A)\cong \HH^*(A,\check A)$
 induced by the isomorphism of $A$ and $\check A$
which is defined by the non-degenerate pairing of $A$.
\end{lem}
\begin{proof}
The only statement to prove is the last isomorphism. As mentioned
the map on the chain level is induced by the isomorphism of $A$ and
$\check A$ defined by the non-degenerate pairing of $A$. The fact
that the complexes are isomorphic follows from the well known fact
that the invariance of the pairing $\la ab,c\ra=\la a,bc\ra$ implies
that the isomorphism between $A$ and $\check A$ is an isomorphism of
$A$ bi--modules, where the bi--module structure of functions $f \in
\check A$ is given by $a'fa''(c)=f(a''ca')$, see e.g.\ \cite{Loday}.
\end{proof}

For any $f\in {\CH}^n(A,A)$ let  $\tilde f\in
\check{A}^{\otimes n}$ be its image under the isomorphism of
$\check A$ with $A$ defined by the Frobenius structure of $A$.

Given pure tensors $f_i=f^{0i}\otimes f_{1i}\otimes \dots\otimes
f_{in_i}
\in {\CH}^{n_i}(A,A), i \in \set{0}{n}$
we write $\tilde f_i= f_{0i}\otimes \dots\otimes
f_{in_i}$ for their image in $\CC^{n_i}(A)$. Fix $\a\in \Anarc(n)$.
 Now decorate the sides $s\in Sides'(\a):=(\amark)^{-1}(1)$
of the complementary regions, which correspond to pieces of the
boundary, by elements of $A$ as follows: for a side $s\in Sides'$
let $j$ its position in its cycle $c_i$ counting only the sides of
$c_i$ in $Sides'$  starting at the side corresponding to the unique
outside angle at the boundary given by the cycle. If the number of
such sides at the boundary $i$ is $n_i+1$  then set $f_{s} :=
f_{ij}$.

\begin{figure}
\epsfxsize = \textwidth \epsfbox{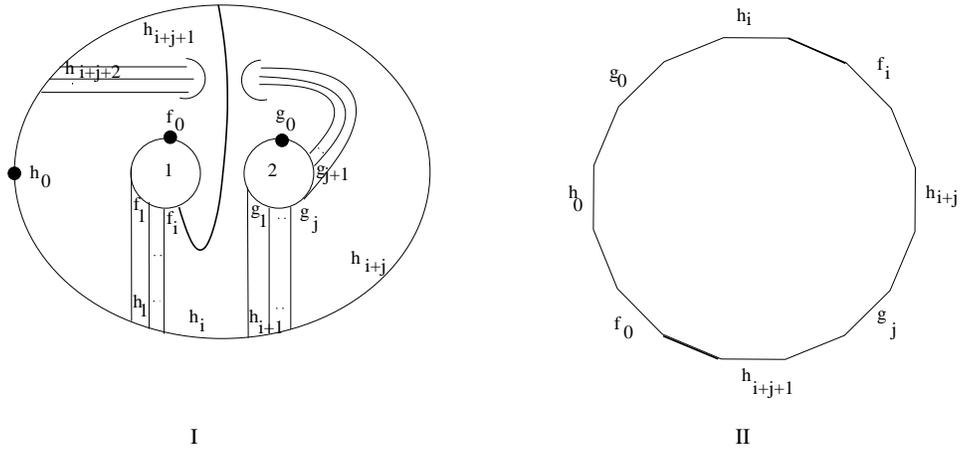}
\caption{\label{g2boundary} A partitioned arc graph with
decorations by elements of $A$ and one of its decorated polygons.
The bold line corresponds to the bold edges.}
\end{figure}

Now we set
\begin{equation}
Y(\a)(f_1,\dots,f_n):=Y(\P(\a))(\bigotimes_{s\in \angle^+(\a^p)}f_s)
\end{equation}
We extend this definition by linearity if $f_i\in
\CH^{n_i}(A,A), i \in \bar n$.  If the condition that
$n_i+1$ equals the number of $Sides'$ at the boundary $i$
is not met, we set $Y(\a)(f_0,\dots,f_n)=0$.
An example of a decorated partitioned surface and its polygons is
given in Figure \ref{g2boundary}.

\subsection{Ribbon correlation functions}
In this section, we give a dual and equivalent picture for the
quasi--filling arc graphs in terms of ribbon graphs. Whether in the
quasi--filling case one wishes to use ribbon graphs or surfaces is
basically a matter of taste. Ribbon graphs seem to be more
established, but actually the surfaces seem to be the better
geometric fit especially if one wants to extend the operations ``to
the boundary'' as explained in the next section. Since CFT is,
however, usually associated with ribbon graphs rather than arc
graphs, we give the details of the construction for this dual
picture. We stress, however, that this duality only exists for
quasi--filling arc graphs indexing cells of $\Anarcn$ and that
furthermore our correlation function are completely different from
those in \cite{Keuro} where basically a CFT is defined from an
$A_{\infty}$--algebra. What we define are morally actually the
correlation functions of the closed string states viewed as
deformations of the category of open strings evaluated on a cell of
the open moduli space, see \S\ref{conclusion}.

\subsubsection{Vertex correlation functions}
Fix an angle marked ribbon graph $(\G,\amark)$. Let $v\in V_{\G}$
be a vertex of $\G$. Let $\angle^+(v)$
be the subsets of angles of $\angle^+(\G)$
whose flags are also incident to  $v$ and
define $Y(v):A^{\otimes \angle^+(v)} \rightarrow k$ by:

\begin{equation}
Y(v)(\bigotimes_{\alpha \in\angle^+(v)} a_{\a}) =\int \prod_{\alpha
\in\angle^+(v)} a_{\a}
\end{equation}
again just as for equation (\ref{picor}) this is well defined as a
function on the tensor product indexed by sets.

\subsubsection{Correlators defined by an angle marked ribbon graph}
Let $\G'$ be a partitioned marked ribbon graph with an angle
marking.

\begin{equation}
\label{ribpartcoreq}
Y(\G')(\bigotimes_{\a\in \angle^+(\G)}a_{\a}) =\bigotimes_{v\in
V_{\Gamma}} Y(v)(\bigotimes_{\beta \in \angle^+(v)}a_{\beta})
\end{equation}

Let $\{c_i\}, i=0,\dots, n-1$
 be the set of cycles labelled by $\set{0}{n-1}$. Let
$\angle^+(c_i)$ be the angles corresponding to the flags of the cycle
$c_i$ which are marked by $1$.
Since the ribbon graph is marked, we have an enumeration of all flags
and hence all the angles,
hence we can  think of the equation (\ref{ribpartcoreq})
as a  map as defined on the subspace
$\bigotimes_{i=1}^n A^{|\angle^+(c_i)|}\subset TA^{\otimes n}$
and extend it by zero outside of this subspace.

Finally for an angle marked, marked ribbon graph, we define
\begin{equation}
Y(\G):=Y(\P(\G))
\end{equation}
by extending linearly.

\subsubsection{Correlators for the Hochschild complex, the ribbon version}
As above for any $f\in \CH^n(A,A)$ let  $\tilde f\in
\check{A}^{\otimes n}$ be its image under the isomorphism of
$\check A$ and $A$ defined by the Frobenius structure of $A$.

Fix  $\G\in \Anrib(n)$, and let $c_i$
be the cycles of the underlying ribbon graph also denoted by $\G$. Set
$n_i=|\angle^+(c_i)|$. Now
for $f_i\in \CH^{n_i-1}(A,A), i \in \set{0}{n}$ which are pure
tensors $\tilde f_i= f_{1i}\otimes \dots\otimes f_{1n_i}$.
Recall that each cycle has a linear order, since the graph was
marked. Now decorate the angles of the graph by elements of $A$ as
follows: for an angle $\a$ let $f$ be the flag of $\a$ and $j$ its
position in its cycle $c_i$ starting at the flag preceding the
marked flag counting only the elements of $\angle^+(\G)$, then set
$f_{\a} = f_{ij}$.

\begin{equation}
Y(\G)(f_1,\dots, f_n):=(\bigotimes_{\a\in \angle^+(\G)}f_{\a})
\end{equation}
We extend this definition by linearity if  $f_i\in
\CH^{n_i-1}(A,A)$ and $n_i=|\angle^+(c_i)|$, otherwise,
if the condition
 is not met, we set $Y(\G)(f_1,\dots,f_n)=0$.

\section{Extension to the boundary and $dg$ properties}
Before we start the discussion of the $dg$-properties, we wish to
point out the following. Let $A$ be a Frobenius algebra and let
$\phi:A \rightarrow \check A$ be the isomorphism defined by the
metric. If $\mu^{\dagger}$ is the adjoint of $\mu$ then  $\phi
\mu^{\dagger}\phi^{-1}=\check \mu:=\Delta$, i.e.\ the natural
coproduct on $\check A$, . Moreover $\phi$ induces an isomorphism
$\psi: TA\to T\check A$. Notice however, this is {\em not an
isomorphism of dg-algebras}, since $\phi \mu \phi^{-1}$ is the
induced multiplication on $\check A$ whereas the natural
differential comes from its co-simplicial structure given by
$\Delta$. We will elaborate on this a little.

As we have discussed for a  Frobenius algebra
 there are canonical isomorphism  $TA\simeq T\check A \simeq \CH^*(A,\check A)
\simeq \CH^*(A,A)\simeq\HC^*(A) $ where we use $\HC$ to indicate
the cyclic cochains. Furthermore the $dg$ structures of the middle two
are compatible if $A$ is a  Frobenius
algebra yielding isomorphisms:
$HH^*(A,A)\simeq HH^*(A,\check A)$
So we can work with $T\check A$ or $TA$ to define the correlators.
In the same vain also the spaces $Hom(TA^{\otimes n},TA^{\otimes m})$ and
$Hom(TA^{\otimes n}\otimes \check TA^{\otimes m},k)$ are isomorphic.
Thus we can also work in the cyclic setting for defining the correlators.

A slight complication arises, when we would like to check the
$dg$--properties of the operadic of PROPic actions defined by
dualizing say $m$ factors of $TA$ as above. The complication is
that although the spaces $Hom(TA^{\otimes n},TA^{\otimes m})$ and
$Hom(TA^{\otimes n}\otimes  TA^{\otimes m},k)$ are isomorphic, if
$A$ is Frobenius, they  have different $dg$--structures when they
are endowed with the natural $Hom$ differentials. In the first
case $\del_{Hom}(f)=f\circ{\del_{TA^{n}}}\mp\del_{T\check
A^m}\circ f$ while in the second case we get $\del_{Hom}(\tilde
f)=f\circ{\del_{TA^{n+m}}}$ where $\tilde f$ the image of $f$
under the isomorphism induced by $\phi$. As discussed above, these
differentials are different. So to get a structure of a $dg$
algebra over a $dg$--PROP, we have to additionally {\em a priori}
declare some boundaries inputs and other outputs. Now on the other
hand, in the geometric models we are considering, say $\Diarc$,
the differential is a topological differential, which is
independent of the discrete structure labelling the boundary. This
independence of the boundary being labelled ``in'' or ``out'' will
be the case for all topological models of surfaces, since these
structures are naturally cyclic and the same type of argument
applies.
 So we
will have to be careful about the type of statements we can make.
We can only expect a compatibility of the $dg$-structure of the
topological chains with the algebraic model if the discrete data
of $\In$ and $\Out$ is canonical. When such a canonical operadically
closed choice of $\In$ and $\Out$ is present, we indeed find the
compatibility of the $dg$ structures.

There is yet another caveat, though. An algebraic complex like the
Hochschild complex does not ``see'' the moduli space structure in
the sense that it does not naturally distinguish any differentials
(say in the co-simplicial setup). In the moduli space case or the case of an
open subset of $\A$ of $\Ana$,
however, certain differentials are set to zero, since we are
dealing with relative chains. In these situations, the Hochschild
differential, will force us to go to the boundary. This type of problem is present for
higher genus and for several ``ins'' and ``outs''.
 In genus zero with only one output it does not appear which explains the naturality of
 the operations of \cite{del,cact} that is in the case of Deligne's conjecture and
 the cyclic Deligne conjecture. In the general case, we have to
grade the  subspaces of  $\CHom(TA)$ and pass to the associated
graded, to obtain the desired operations.

\subsection{$Dg$ properties of the PROP action}
\subsubsection{The tree level operads}
If we fix one output boundary and chose the natural embedding of
$\Tree_{cp}$ into $\Diioarc$,
 then we are in the case of \cite{del,cyclic} and the
action is indeed the $dg$ action discussed for $\Lintree_{cp}$ and
$\Tree_{cp}$, respectively\footnote{see \cite{hoch1} for the
definitions of these subsets.}.

\subsubsection{The Sullivan-Chord diagram case}
In the case of $\Gr\OC(\ioarc)$ with the conventions of
\cite{hoch1} reviewed in the Introduction, we can again put
ourselves into  the setting of a $dg$-action of a $dg$--PROP on a
$dg$--algebra, but we need to extend the action  to the boundary.
Generalizing the arguments of \cite{del,cyclic} we will show that
the differential on the Hochschild side corresponds to the
differential of $\A$ restricted to $\Diioarc$. If one removes an
arc from  a graph indexing a cell in $\Gr\OC(\ioarc)$, it is not true in general, that we still
obtain a graph indexing a cell of $\Gr\OC(\ioarc)$. This does not
hold even if the genus is zero, so we will have to extend the
action to the graded cells of $\Diioarc$. Fortunately there is a
CW-complex which models which naturally allows us to do this, that is  $CC_*(\Diioarci)\cong
\Gr\OC(\Diioarc)$.

We will see that in order to achieve a well defined action, we
will need some additional assumptions. These are satisfied if we
restrict our attention to  a  {\em commutative} Frobenius algebra
$A$.

{\sc Assumption:} For the rest of the
discussion of this subsection let $A$ be a commutative Frobenius algebra.

\subsubsection{Extending to the boundary}
\begin{figure}
\epsfxsize = \textwidth \epsfbox{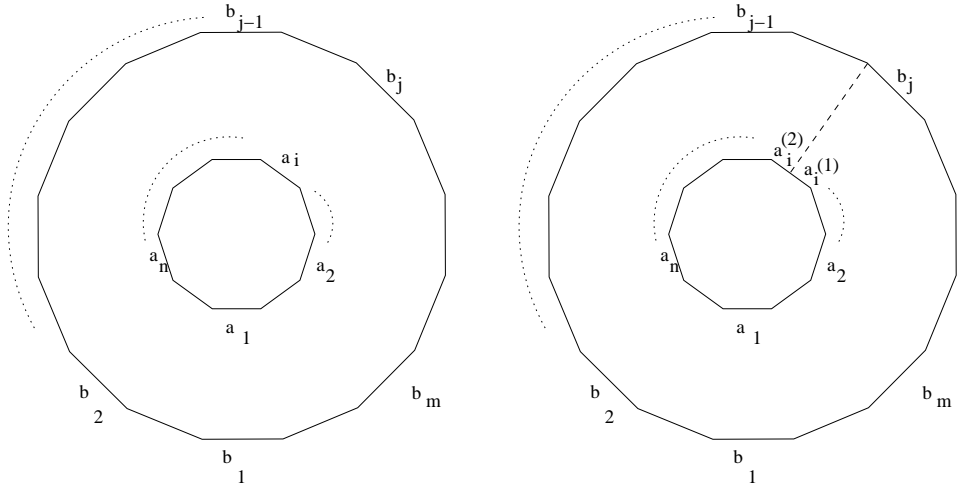}
\caption{\label{gluecyl} A cylinder with boundary components
labelled by elements of $A$ and an indicated cut.}
\end{figure}

The definition of this extension is dictated by the $dg$
condition. This means that we will have to consider elements in
$\Arc$ which are in the limit of elements of $\Arcno$. Removing an
arc corresponds to gluing together two polygons, and so we have to
deal with not only polygons, but also with cylinders and so forth.
For a cylinder $C(n,m)$ with two boundaries given by polygons
$p_1\in \Poly(n)$ and $p_2\in \Poly(m)$,  and a choice of cut
indices $(i,j)$, see Figure \ref{gluecyl} for an example, we
define

\begin{multline}
\label{ceq} Y(C(n,m),(i,j))(a_1,\dots, a_n,b_1,\dots,b_m):=\\
\int\sum  a_1
\dots a_i^{(1)} b_j \dots b_m b_1 \dots b_{j-1} a_i^{(2)} \dots
a_n
\end{multline}
where $\Delta(a_i)=\sum a^{(1)}\otimes a^{(2)}$ is the co-product
of $a_i$ using Sweedler's notation.

This definition is forced on us, if we wish to ensure compatibility with the Hochschild
differentials.

 \begin{rmk} In order to ensure that equation (\ref{ceq}) is
independent of the choice of cutting edge that is the
indices $i$ and $j$ in equation (\ref{ceq}) we assumed commutativity.
\end{rmk}

\begin{lem}
\label{transfer} If $A$ is a commutative Frobenius algebra then
the r.h.s. of equation (\ref{ceq}) is independent of the choice of
$i$ and $j$ and coincides with
\begin{equation}
Y(n,m)= \int a_1 \dots a_n b_1\dots b_m e
\end{equation}
 where $e=\mu \circ
\Delta(1)$ is the Euler element.
\end{lem}

\begin{proof}
If $A$ is commutative, then we have
$$
Y(C(n,m))( a_1,\dots, a_n,b_1,\dots,b_n):=\int  a_1 \dots a_{i-1}(\mu
\circ \Delta) (a_i)a_{i+1} \dots a_n b_1 \dots b_m
$$
but in any commutative Frobenius algebra one has
\begin{multline*}
\int (\mu \circ \Delta) (a)bc=\la(\mu \circ \Delta) (a),bc\ra\\ = \la
a,(\mu \circ \Delta) (bc)\ra =\la a,b(\mu \circ \Delta )(c)\ra =
\int ab(\mu \circ \Delta) (c)
\end{multline*}
fixing $a=a_i, b=\prod_{k\neq i}a_i\prod b_j$ and $c=1$ shows the claim.
\end{proof}

\subsection{Correlators for $\Ana$}
In general we extend the action as follows. Notice that given an arc
graph $\a$ each complementary region $S\in \comp(G)$ has the
following structure: it is a surface of some genus $g$ with $r\geq1$
boundary components  whose boundaries are identified with a 2k-gons.
Alternating sides belong to arcs and boundaries as above and the
sides come marked with $1$ or $0$ by identifying them with the
angles of the underling arc graph. Now let $Sides'(S)$  be the sides
which have an angle marking by $1$ and let $\chi$ be the Euler
characteristic of $S$.  We set
\begin{equation}
\label{surfacecor} Y(S)(\bigotimes_{s\in Sides'(S)}a):=\int
(\prod_{s\in Sides'(S)} a_s) e^{-\chi+1}
\end{equation}
where $e:=\mu(\Delta(1))$ is the Euler element. For an angle
marked partitioned arc graph $\a^p$ we set

\begin{equation}
\label{partcor} Y(\alpha^p)(\bigotimes_{S\in \comp(\alpha_i)}
(\bigotimes_{s\in Sides'(S)}
a_s))=\bigotimes_{S\in \comp(\alpha_i)}Y(S)(\bigotimes_{s\in
Sides'(S)} a_s)
\end{equation}

Again, for $\a\in CC_*(\Ana)$ we simply set
\begin{equation}
\label{boundaryydef}
Y(\a)=Y(\P(\a)).
\end{equation}

\subsubsection{The Hochschild differential}
Consider the operation $Y(\a^p)$ for $\a\in \Gr\OC(\ioarc_{\#})$. The $Hom$
differential on this viewed as an element in
$Hom(CH^{\In}(A,\check A),CH^{\Out}(A,\check A)$
is given by $\del_{Hom}(Y)(f_i)=Y(\del_{Cyc}(f_i))\mp \del_{Cyc}Y(f_i)$.
Here we indexed the tensor products by the sets $\In$ and $\Out$ and denoted
the differential of the cyclic bar complex by $\del_{Cyc}$.

We can consider $Y\in Hom(CH^{\In}(A,\check A)\otimes
(CH^{\Out}(A,\check A)^{\vee},k) $, by decorating the $\In$
boundaries with the elements $f_i$ and the $\Out$ boundaries by
elements $a_i\in A$. That is $Y(f_0,\dots, f_n)(\bigotimes a_i)$.
Then the first term in the differential is given by applying
$\Delta=\check \mu$ cyclically to each element $f_i$ in the left
hand side viewed as an element in $T\check{A}$ decorating the
$\In$ boundaries. The second term in the differential is given by
the sum obtained by decorating exactly one of the angles of the
$\Out$ boundaries with the product of two variables $a_ia_{i+1}$.
These summands will cancel with summands from the first term
essentially due to the Frobenius condition
\begin{equation}
\label{parallelcancel}
\la f_{ij},a_k a_{k+1}\ra=\la \Delta(f_i^j),a_k \otimes a_{k+1}\ra
\end{equation}
where we wrote $f_{ij}$ using the notation of \ref{hochactionpar}.
Hence,
we are left with the summands of the first term that are not cancelled.
These correspond to
angles on the $\Out$ boundaries marked by $0$.

\begin{figure}
\epsfxsize = \textwidth
\epsfbox{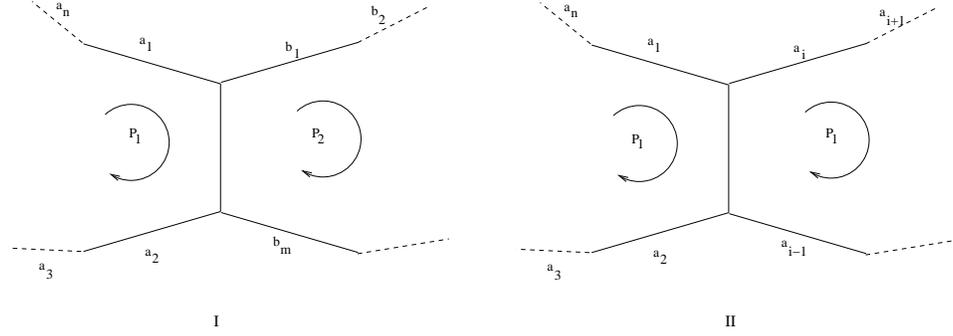}
\caption{\label{sepnonsep}
I. The separating case. II. The non--separating case}
\end{figure}

More precisely consider decorating two neighboring angles at an
$\In$ boundary by say $\Delta(f_{ij})$ where the common edge $e$
belongs to an angle marked by $0$ on the $\Out$ boundary. There
are two cases. Either the edge is separating, that is it separates
two different complementary regions or it is non-separating, that
is the same region lies on both sides of the edge, see also Figure
\ref{sepnonsep}.

First let's consider the angle markings all given by $1$.
Let $P_1$ and $P_2$ be the two complementary regions
on the two sides
of the edge $e$.
And fix the notation $a_1,\dots, a_n$
 for the elements decorating the sides of the polygon
$P_1$ and $b_1,\dots, b_m$ for the elements decorating
 the sides of $P_2$ where in both cases the enumeration
 is compatible with the cyclic order. Also
 set $\Delta(f_i^j)=\sum a^{(1)}\otimes b^{(2)}$
 in Sweedler's notation. Say  that $a^{(1)}=a_1$ and $b_1=b^{(2)}$;
 see Figure \ref{sepnonsep}.
 In the separating case applying the
Hochschild differential before applying $Y$,
we obtain a contribution of the type
\begin{multline}
\label{sepcor}
\sum \int  a^{(1)} a_2  \dots a_n \cdot \int b^{(2)} b_{2}\dots b_m
= \int f_i^j b_2 \dots b_{m-1} (b_m a_2) a_3\dots a_n
\end{multline}
which is the contribution obtained by gluing the two polygons along the common edge and decorating one of the two joined sides by $f_i^j$  and the other by $b_m a_2$.

In the non-separating case, we only have one complementary region
$P_1$, let $a_1, \dots, a_m$ be the elements decorating the sides.
Also again use $\Delta(f_i^j)=\sum a^{(1)}\otimes a^{(2)}$  and
let $a_1=a^{(1)}$ and say $a_i:=a^{(2)}$ the contribution reads

\begin{multline}
\label{nonsepcor} \pm \int  a^{(1)} a_2 \dots
a_{i-1}a^{(2)}_ia_{i+1} \dots a_n= \int  (\mu\circ \Delta) (f_i^j)
a_2 \dots a_{i-1} a_{i+1} \dots a_n\\= \int  f_i^j a_3 \dots
a_{i-1}a_{i+1} \dots (a_na_2) e
\end{multline}
where $e=\mu\circ\Delta(1)$.

First if both polygons are quadrangles  with all markings $1$ then
the term in equation (\ref{sepcor}) cancels with the Hochschild
differential on $Y$ as discussed above. In all other cases, if
$\a$ is an quasi-filling element of $\Diioarc$ endowed with the
standard angle markings then one of the elements $b_m,a_2$ is one
in equation (\ref{sepcor}), so that we indeed obtain the
contribution to the correlation function associated to the surface
in which the arc has been removed. In the same situation but with
a separating arc in equation (\ref{nonsepcor}) either $a_2$ or
$a_m$ is equal to one and we again obtain the contribution to the
correlation function associated to the surface in which the arc
has been removed --- now basically by definition.

This calculation generalizes  to arbitrary $Y(\a^p)$ for $\a\in
\Gr\OC(\Diioarc)$.

\begin{prop}
The $dg$--structures of $CC_*(\Diioarci)$ and $\Hom(\CH(A,A)$ are
respected by $Y$. Or in other words: The equation
(\ref{surfacecor}) and hence the equation (\ref{partcor}) are the
operations corresponding to the summands of the $\Hom$
differential of $Y(\a^p)$ and furthermore these summands
correspond to the respective boundary components of $\a$.
\end{prop}

\begin{proof}
In the case of splitting angles, we see by
the considerations above
that there are two terms in the $Hom$ differential which cancel.
In fact all of the terms $\del \circ Y$ cancel in this way.
On the other hand
generalizing the formalism explained above, we see that when two
neighboring angles which are not parallel are assigned the
coproduct of an element $a_i$, the resulting operation is the
operation associated to the arc graph in which the arc
corresponding to the common edge is removed. This either causes
 two bordering complementary regions to be joined or self--glues a
 complementary region to itself. In both cases
  the resulting function is the product of integrals over all the
boundaries of the joined surface. In the case of self--gluing this
yields a term $\mu(\Delta(a_i))$. As in Lemma \ref{transfer}, we
can ``transfer'' the $\mu\circ \Delta$ to an inserted unit. Now
iterating this process, we have to remove $-\chi+1$ edges to
obtain a complementary region $S$ with genus $g$ and $r$ boundary
components. This accounts for the tensor power of $e$. Also
iterating the argument of the ``transferring'' the
$\mu\circ\Delta$  from the elements $a_i$ to $1$ one obtains the
independence of the exact incidences of the removed edges. I.e.\
if the surface $S$ can obtained by removing other edges, the
resulting operation will be the same. This means that such
an expression in the differential is well defined.
\end{proof}

Collecting the results, we find:

\begin{thm}
The $Y(\a)$ defined in equation (\ref{boundaryydef}) give operadic
correlation functions for $CC_*(\Diioarci)$ and induce a
$dg$--action of  the $dg$-PROP $CC_*(\Diioarci)$ on the
$dg$--algebra $\overline {CH}^*(A,A)$ of reduced Hochschild
co-chains for a commutative Frobenius algebra $A$.

The $Y(\a)$ also yield correlation functions on the tensor algebra
of the co-cycles of a differential algebra $(A,d)$ over $k$
 with a cyclically invariant trace $\int: A\to k$
that satisfies $\int da=0$ and whose induced pairing on $H=H(A,d)$
turns $H$ into a Frobenius algebra. These correlations functions are
operadic chain level correlation functions.
\end{thm}

\begin{proof}
We have shown in \cite{hoch1} that the $CC_*(\Diioarci)$ form a
$dg$--PROP that is isomorphic to the $dg$--PROP $\Gr\OC(\Diioarc)$
and we have that the map $\PA$ is operadic/PROPic. Now the gluing
for the correlation functions coincides with the algebraic one
used to define the PROP structure of Proposition
\ref{anglegluingopprop}. On the Hochschild side, we simply plug in
elements. On the graph side, we correspondingly plug in angles
marked with $1$ to angles marked with $1$. The first thing, we
have to make sure is that the steps (1) and (2) in Definition
\ref{partgluedef} are respected on the Hochschild side. By
\cite{cact} the double twisted case of step (1) corresponds to
applying the Connes' operator $B$ to both sides, that is applying
$B^2$,
 and hence yields zero.  The case of (2) cannot occur since the arcs
are only running from $\In$ to $\Out$. Therefore the gluing
actually corresponds to the ``twisted'' gluing. Lastly, in the gluing
for $\PA\Diioarc$ there are no terms of lower degree, since
one never glues separating to separating angles, since these
are labelled by $0$ on the $\Out$ boundaries and by $1$ on the $\In$
boundaries. So the gluing on the Hochschild side corresponds
to the gluing in the associated graded for the non--partitioned graphs
before applying $\PA$.

The last statement follows from Proposition \ref{dgcorrelatorprop}
\end{proof}

\begin{rmk}
Seemingly related results have been obtained by \cite{TZ} in a different setting.
Their definition
of Sullivan Chord diagrams is, however, different from
ours and, as far as we can see, also the  definition
of the action also differs. It is therefore not possible to relate their calculations
to the present ones or those of \cite{cyclic,del} directly. It
would be interesting to know how if it is possible to compare the
two actions despite their different settings.
\end{rmk}

\begin{cor}
The operadic correlation functions descend to give a PROP
action of $H_*(\Diioarci)$ on $\HH^*(A)$ for a commutative Frobenius
algebra $A$.
\end{cor}

\subsection{Co-simplicial properties of the action of moduli space}

\subsubsection{The operation of $\LDiioarc$}
We recall from \cite{hoch1} that $\LDiioarc$ is the subspace
of $\Diioarc$ whose underlying arc graphs are not twisted at the $\In$
boundaries. Using the constant marking $\amark \equiv 1$ this
space is a subspace of $\Ana$.

\begin{prop}
The correlation functions  (\ref{boundaryydef}) are operadic
correlation functions for the PROP $\OC(\LDiioarc)$,
 the tensor algebra on the co--cycles $Z(A)$ of a quasi--Frobenius
algebra. That is they give chain level correlators.
\end{prop}

\begin{proof}
The algebraic operadic composition on the level of partitioned
angle marked arc graphs corresponds to the insertion on the
$\CHom(TZ(A))$ side. The conditions for the twisted gluing needed
to make $\PA$ operadic never occur. There are never any double
twisting and also never any closed loops.
\end{proof}

\subsubsection{The tree level: $\Tree$}

As we have previously discussed, one cannot expect that the cyclic
operad will go over to the $dg$ setting. We, however, have the
following interesting observation.

\begin{lem}
\label{mccompat} Using the isomorphism $\CH(A,A)\simeq \overline{TA}$
for a Frobenius algebra,
the operations defined by the cells of $\Tree$ embedded into
$\Anarcn$ by the constant marking $\amark\equiv 1$ yields the
operations $\sqcup$ and $\square$ of \cite{MScosimp}. In
particular,  the operation of $\LDiioarc$ induces the $\Xi_2$
operation of \cite{MScosimp}.
\end{lem}

\begin{figure}
\epsfxsize = \textwidth
\epsfbox{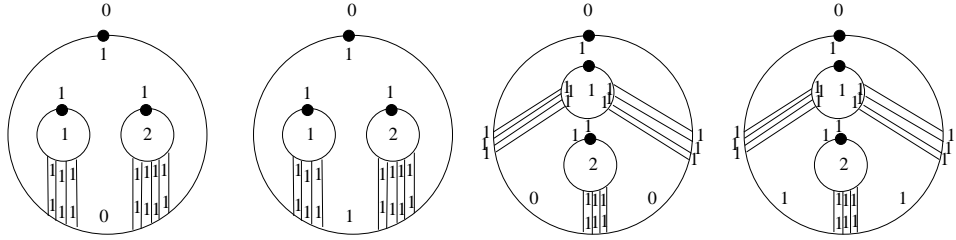} \caption{\label{partitioned}
Examples of the angle marked partitioned families yielding $\cup$,
$\sqcup$, $\circ_i$ and
$\square_i$}
\end{figure}

\begin{proof}
This is a straightforward verification. The relevant arc families
are depicted in Figure \ref{partitioned}.
\end{proof}

\begin{rmk}
In other words using the language of \cite{MScosimp}, our operations
allow us to recover the operation of the functor operad given by the
Kan extension of the operad of sequences with differentials up to
complexity 2. This action corresponds to the genus--$0$--$\Lintree$
operad embedded by the marking $\amark \equiv 1$. In general, even
going to the boundary as in the last paragraph, we cannot expect to
get actions of the higher $\Xi_n$. This is commensurate with
Deligne's conjecture. Since in the differential structure of the
sequences of higher complexity, one does not retain the topological
information of the underlying surface when moving to the boundary,
one would, in contrast to the last paragraph, have to identify the
surfaces $S$ again with polygons, which is not true in $\Arc$, but
can be done in $\mathcal{S}t\mathcal{A}rc$, the stabilization of the
arc operad.

 Two questions arise in this setting. What are the
conditions to get the higher differentials and is there a
co-simplicial type of setup for $\Arcno$ or $\Anarcno$. The first
question will be addressed in \cite{Ribbon} where we will deal
with the stabilization of the arc operad and its loop structure.
\end{rmk}

\subsection{New Monoidal structures on the Hochschild co-chains}
In order to match the geometric grading by non-partitioning angles
on the algebraic side, we will introduce a new monoidal structure
of $A$--bi--modules, where $A$ is an associative ring. This structure
will allow us to define a natural grading on a subcomplex
of the Hochschild--complex whose associated graded is the correct
receptacle of our correlation functions.

\begin{df}
Let $A$ be an associative ring and let $M$ and $N$ be
$A$--bi--modules. We define $M\boxtimes N$ to be $M\otimes
A\otimes N$ with the obvious structure of $A$--bi--module.
\end{df}

\subsubsection{A new co-product}
Now if we study the basic operations on the Hochschild co-chains,
which are a generalization
of the operations of $\Shuff$ to the brace sub--operad $\Brace$,
 we see that they naturally correspond to slightly different
operations induced by the new monoidal structure.

We will first treat the operad corresponding to $\Shuff$. For this
we  notice that with the monoidal structure $\boxtimes$  on $TA$
the product is a map $\boxtimes: TA\boxtimes TA\to TA$ and dually
the co-product is a map $\lozenge: TA\mapsto TA\boxtimes TA$. Here
the first map is given by $(a_1\otimes \dots \otimes a_n)\otimes b
\otimes (c_1\otimes \dots \otimes c_m) \mapsto a_1\otimes \dots
\otimes a_n\otimes b \otimes c_1\otimes \dots \otimes c_m$, that
is it ``raises degree by one''; and the second map is given by
$a_1\otimes \dots \otimes a_n \mapsto \sum_i \pm (a_1\otimes
a_{i-1})\otimes a_i \otimes (a_{i+1}\otimes \dots \otimes a_n)$.

Now the multiplication $\cup$ in $CH^*(A,A)$ was given by the
multiplication $\otimes$ in the tensor algebra $T\check A$, and
the multiplication $\mu$ of the algebra $A$.

\begin{equation}\cup:  A\otimes \check A^{\otimes n}\otimes A\otimes
\check A^{\otimes m}\stackrel {\otimes \circ\sigma }{\rightarrow}
A\otimes A \otimes \check A^{\otimes n+m} \stackrel{\mu
}{\rightarrow} A \otimes \check A^{\otimes n+m}
\end{equation}
here $\sigma$ is just the permutation of the tensor factors and
$f\cup g(a_1,\dots a_n,b_1,\dots b_m)=f(a_1,\dots,
a_n)g(b_1,\dots, b_m)$.

Using the new monoidal structure for the same canonical maps we
obtain a new multiplication

\begin{equation}
\sqcup: (A\otimes \check A^{\otimes n})\boxtimes (A\otimes \check
A^{\otimes m})\stackrel {\boxtimes \circ\sigma }{\rightarrow}
A\otimes A \otimes (\check A^{\otimes n}\boxtimes \check
A^{\otimes m}) \stackrel{\mu }{\rightarrow} A \otimes \check
A^{\otimes n+1+m}
\end{equation}
here $\sigma$ is again just the permutation of the tensor factors
and
\begin{equation}
f\sqcup g(a_1,\dots a_n,b,c_1,\dots c_m)=f(a_1,\dots,
a_n)bg(c_1,\dots, c_m)\end{equation}

This is exactly the operation induced by
the co-simplicial structure used in \cite{MScosimp}, which we
recover using the embedding $\amark \equiv 1$.

In the setting of operadic correlation functions, we are using the
coproduct to separate the different tensor factors of $\check A$
at the different boundaries before integrating over them.

This means that in the current setting, we should  again use the
co-product $\lozenge$. Now each time we use the new co-product
$\lozenge$ this has the effect of inserting a tensor factor of
$A$. So that for instance the usual $\circ_i$ operations of
$CH^*(A,A)$ which use two co-products become operations
$\square_i$ where
\begin{multline}
\square_i(f,g)(a_1,\dots,
a_{n+m+2})=\\
f(a_1,\dots, a_{i-1},a_ig(a_{i+1},\dots
a_{i_m})a_{i+m+1},a_{i+m+2},\dots,a_{n+m+2}).
\end{multline}

 This is again the
operation obtained by \cite{MScosimp}.

\subsection{Graded correlators}
Now in general from the correlators of $\a^p$ whose underlying
arc graph $\a$ is in $\OC(\Arcn(F))$  we obtain maps
\begin{equation}
Y(\a^p): TA^{\otimes n}\to \bigotimes_{i=1}^n A^{\boxtimes n_i+1} \to k
\end{equation}

\begin{nota}
\label{lozengenota} If $k=|E(\a)|$, and $p\in P(n,k)$, we will use
the following notation: $\lozenge^{l}:TA\to TA^{\boxtimes l+1}$ is
the iteration of $\lozenge$ given by $(\lozenge \otimes
(id_A\otimes id_{TA})^{\otimes l})\circ \lozenge \otimes
(id_A\otimes id_{TA})^{\otimes l-1}) \circ \dots \circ (\lozenge
\otimes id_A\otimes id_{TA}) \circ \lozenge$.
\end{nota}

\begin{equation}
\label{corrdecomp}
Y(a^p)= (\bigotimes_{\pi \in Comp(\a^p)} Y_{Poly_2}(\pi)
\otimes \bigotimes_{i=1}^k\eta^{n_i-1})\circ \sigma \circ
\bigotimes_{i=0}^n \lozenge^{n_i}
\end{equation}
where the we think of the complementary regions $Comp(\a)$ as a
subset of the complementary regions of $\a^p$ and $Y_{Poly_2}$ are
the polygon correlation functions defined in equation
\ref{polycoreq}. Here $\sigma$ permutes the factors of $TA$ and
the factors of $A$ corresponding to $\bigotimes_{i=1}^n
TA^{\boxtimes n_i+1}$ and we used Notation \ref{etanota} and
Notation \ref{lozengenota}.

\begin{df}
We let $\Modshuffint(n)\subset \Hom(TA)(n)$
be the image of all the operations of
$\PA\Arcn(n)$, by considering $0$ as ``out''.
\end{df}

We get an analogous statement to the Proposition \ref{modshuffprop}.

\begin{prop}
\label{modshuffintprop} After dualizing to obtain elements in
$Hom(TA^{\otimes n+1},k)$ any element in $\Modshuffint(n)$ can be
written uniquely as in equation (\ref{corrdecomp}). Set
$l=\frac{1}{2}\sum_i (n_i+1)-1$ then $\Modshuffint$ is graded by
$l$. Moreover the composition in $\Modshuffint$ respects the
induced filtration of elements of degree $\leq l$. Lastly,
 the decomposition identifies $\Modshuffint$
with the subspace of $\Hom(TA)$ obtained by dualization
for the subspace generated by the coproduct $\lozenge$,
permutations of the factors $A$ and $TA$, and $\eta$
in $\bigoplus_n Hom(TA^{\otimes n+1},k)$.
\end{prop}

\begin{proof} Completely analogous to the proof of \ref{modshuffprop}.
The first statement is again clear by the definition of
$\Modshuffint$ as the image. Likewise, the last statement is also
again straightforward, by arranging the operation in the specified
order. On the other hand it is easy to give the arc graph in
$\P\Anarcno(n)$ by drawing one arc for each factor of $\eta$ with
the incidence relations given by $\sigma$. One quick way is to
dually draw the partitioned ribbon graphs which one vertex per
complementary region. This identifies the two subspaces. In this
identification there is one factor of $\lozenge$ for each inner
angle which is not partitioning. The last claim, that the
operations respect the filtrations is clear after identifying $k$
with the dimension of the cell, that is the number of edges minus
one, of the underlying graph for the operation. The mentioned
equality follows from the combinatorial identity $|\angle_{\rm
inner}|+|\angle{\rm outer}|=|{\rm Flags}| =2 |{\rm edges}|$ which
still holds true.
\end{proof}

\begin{prop}
\label{tamodshuffintprop}
 For any Frobenius algebra $A$
 the correlations functions of equation (\ref{corrdecomp})
 define operadic correlation functions for $\Gr\PA\Arcn$
with values in  the associated graded $\Gr\Modshuffint$ of
$\Modshuffint\subset\Hom(TA)$.
 By regarding $\a\to Y(\PA\a)$ the same statement hold
also for   $\Gr\OC(\Anarc)$.
\end{prop}

\begin{proof}
Using the Casimir element to dualize,  we only have to show that
the resulting structure is that of an algebra over a cyclic
operad. Again the $\Snn$ equivariance is manifest. After dualizing
the gluing on the flags in the operadic composition turns into the
identity map $id: A\to A$, so that indeed the gluing $\circ_i$ on
$\PA\Arcn$ maps to insertion at the $i$-th place in $\Hom(TA)$.
Dealing with the extra steps (1) and (2) in the definition of the
gluing in $\PA\Arcn$, we see that on the side of $\Hom(TA)$ they
would not yield zero. However in $\Gr\Modshuffint$ both these
cases are projected out, since they correspond to operations of
lower degree. The same is true for both $\Gr\PA\Anarcno$ and
$\Gr\OC\Anarcno$. by definition.
\end{proof}

Defining the action on the Hochschild complex trough the tensor
algebra as in \S\ref{hochactionpar}, we obtain:

\begin{thm}
Let $A$ be a Frobenius algebra and let $\CH(A,A)$ be the
Hochschild complex of the Frobenius algebra, then the cyclic chain
operad of the open cells of $\Anarc$ act on $\CH(A,A)$ via
correlation functions. Hence so do all the suboperads,
sub-dioperads and PROPs of \cite{hoch1} mentioned in the
introduction. In particular the graph complex of $\Mngn$, the
Moduli space of pointed curves with fixed tangent vectors at each
point act on $CH(A,A)$ by its two embeddings into $\Anarcno$.
Furthermore, on $\PA\Arcno$ the correlation functions are operadic
correlation functions with values in $\Gr\Modshuff$. Moreover, the
operations of the suboperad $\Tree_{cp}$ correspond to the
operations $\sqcup$ and $\square_i$ induced by $\Xi_2$ as defined
in \cite{MScosimp}.

The same formula equation (\ref{corrdecomp}) also yields  operadic
correlation functions for the tensor algebra of the co-cycles of a
differential algebra $(A,d)$ over $k$ with a cyclically invariant
trace $\int: A\to k$ which satisfies $\int da=0$ and whose induced
pairing on $H=H(A,d)$ turns $H$ into a Frobenius algebra, i.e.\ they
are chain level operadic correlation functions with values in
$\Gr\Modshuff$.
\end{thm}

\begin{proof}
We use operadic correlation   function $Y$ above  to give maps
$CH^{p_1}\otimes \dots \otimes CH^{p_{n+1}}\simeq \check A^{\otimes
p_1+1} \otimes \dots \otimes \check A^{ p_n+1} \rightarrow k$. All
the necessary properties follow from Proposition \ref{mccompat} and
Proposition \ref{tamodshuffintprop}. The last statement again
follows from Proposition \ref{dgcorrelatorprop}.
\end{proof}

\subsection{Application to String-topology}

Let $M$ be a simply connected compact manifold $M$ and denote  the
free loop space by $\mathcal{L}M$ and let $C_*(M)$ and $ C^*(M)$
be the singular chains and (co)-chains of $M$. We know from
\cite{jones,CJ} that $C_*(\mathcal{L}M)=\CH^*(C^*(M,C_*(M))$ and
$H_*(\mathcal{L}M)\simeq \HH^*(C^*(M),C_*(M))$. Moreover $C^*(M)$
is an associative $dg$ algebra with unit, differential $d$ and an
integral ($M$ was taken to be a compact manifold) $\int:C^*(M)\to
k$ such that $\int d\omega=0$. By using the spectral sequence and
taking field coefficients we obtain operadic correlation functions
$Y$ for $\Tree$ on $E^1=\CH^*(H,H)$ which converges to
$\HH^*(C^*(M))$ and which induces an operadic action on the level
of (co)-homology. Except for the last remark, this was established
in \cite{cyclic}.

\begin{thm}When taking field coefficients,
the above action gives a $dg$ action of a $dg$--PROP of
 Sullivan Chord diagrams on the $E^1$--term of a spectral sequence
converging to $H_*(LM)$, that is the homology of the loop space a
simply connected compact manifold and hence induces operations on
this loop space.
\end{thm}
\begin{proof}
Recall from \cite{CJ} the isomorphism
$C_*(\mathcal{L}M)=\CH^*(C^*(M,C_*(M))$ comes from dualizing the
isomorphism $C_*(\mathcal{L}M)=\CH_*(C^*(M))$\cite{jones}.
Calculating the latter with the usual bi-complex \cite{Loday} then
we see that the $E^1$-term is given by $CH_*(H^*(M))$ and
dualizing the corresponding $E^1$ spectral sequence, we get
$CH^*(H^*(M),H_*(M))$, so we get an operation of the $E^1$ level.
Since the operation of $\Tree$ was $dg$, it is compatible with the
$E^1$ differential and hence gives an action on the convergent
spectral sequence computing $H_*(\mathcal{L}M)$ and hence on its
abutment.
\end{proof}

\section{Concluding remarks}
\label{conclusion} In this paper and its first part \cite{hoch1}
we have systematically used the $\Arc$ operad and its cousins to
give operations on the Hochschild co-chains of a Frobenius
algebra, by extending and building on our results of
\cite{del,cyclic}. In particular, we have given correlation
functions for $\Arcn$. In physics terms this could be expected by
using the logic of \cite{KR} as follows. If the closed string
states are thought of as deformations of the open string states
and the open string states are represented by a category of
$D$-branes, then the closed strings should be elements of the
Hochschild co--chains of the endomorphism algebra of this
category. Now thinking on the worldsheet, we can insert closed
string states. That is for a world sheet, we should get a
correlator by inserting, say $n$ closed string states. This is
what we have done, if one simplifies to a space filling $D$-brane
and twists to a TCFT.

For string topology, we have given operations using the spectral sequence,
so the question remains, if we loose any information
by passing to the associated graded. This is indeed an interesting question.
It seems though that since all the operations of string topology
 preserve the grading and not just the filtrations, we have
not lost any information.
A question that one could ask is how different possible lifts
 from the associated graded to the filtered complex are related.
It is conceivable that an interesting ``up to homotopy'' structure
is lurking which may possibly be related to Frobenius manifold
structures found e.g.\ in \cite{Merkfrob}. An interesting
observation in this respect is that the operadic correlation
functions allow one to lift to at least the co--cycle level in the
tensor algebra setting. Perhaps this gives enough information to
compare the two sets of operations. It seems that comparing to
\cite{CJ} the operations should even be the same. Although {\it a
priori} they might differ, the operation of \cite{CJ} do not only
respect the filtration, but they act with a definite bi--degree in
the bi--grading and hence {\it a posteriori} seem to have no lower
order contributions in the filtration.

It seems that the combinatorial version of the moduli space of
\cite{P2,KLP} is particularly suited for these applications. One
amazing coincidence is that the Hochschild differential forces one
to consider Penner's compactification. That is it forces to move
from $\ioarc$ to $\Diioarc$. Another interesting remark is that
the grading by the number of arcs on the Hochschild side is
reminiscent of the operation of open strings rather than closed
strings in the framework set up in \cite{KP}. Here we have the
additional restriction that there is exactly one arc per window in
the terminology of \cite{KP}, and the $D$-brane label corresponds
to dualizing the respective element labelling the marked point to
live in $A$. This observation could be a ``shadow'' of the
open/closed duality.

The interplay between algebra and geometry is astonishing, the
algebra side for instance demands the insertion of degeneracies in
order to obtain the $BV$--operator $B$. This manifests itself in
the restriction to the PROP $\ioarc$ and the preservation of the
 algebraic $dg$--structure on the geometric side then forces one to move
to the boundary, viz.\ $\Diioarc$ or the Sullivan--Chord diagrams.

This leads us to interesting aspect which we have left untreated
is the co--simplicial setup for the $\Arc$ complex. That is
reverse engineering the $\Arc$--complex, by starting with a
cosimplicial model coming from partitioned arc graphs. For the
subset of $\Lintree$ this is essentially what has been done in
\cite{MScosimp}. Hence one could expect that the totalization of
the arising complex operates on the totalization of the relevant
Hochschild complex. If this is possible one would have the hope of
``$dg$ compatibility'' after passing to the totalization. One of
the difficulties, however, is that there is more than one
``output'', so that one cannot directly use a co--simplicial
structure since this relies on the category of maps, viz.\ several
inputs, but only one output. In order to accommodate this one
either has to break the cyclic setting of the cyclic operad or one
has to construct a suitable category of sets with correspondences.

Of course a generalization to the $\Ainf$ case would be very useful.
The example of polygon correlators for $\Ainf$--algebras shows a possible path.
The tree level version will be worked out along these lines in \cite{KSch}.
A further area which deserves  study
are the implications for the operations on the cyclic co--chains
and the associated $S^1$ equivariant theories, e.g.\ in the spirit of
\cite{craig}.

Lastly, we wish to point out that at several points we had to
avoid closed loops. On the topological side this basically comes
from the cell decomposition of moduli space, which does not have
any graphs with closed loops. On the Hochschild side this was not
as natural. We avoided the occurrence, by either restricting the
type of graph or passing to the associated graded. In the setting
of partially measured foliations however these closed leaves are
very natural. So one cannot help but wonder if there is yet
another generalization of this whole story to foliations as
outlined in the Appendix of \cite{KP}.

\end{document}